\documentclass{amsart}
\usepackage{amssymb}
\usepackage{graphics}

\theoremstyle{plain}
\newtheorem{lemma}{Lemma}

\newtheorem{proposition}{Proposition}
\newtheorem{theorem}{Theorem}
\theoremstyle{definition}
\newtheorem{definition}{Definition}
\newtheorem{corollary}{Corollary}
\theoremstyle{definition}
\newtheorem{remark}{Remark}

\begin{document}
\title[Singular unitarity in ``quantization commutes with reduction'']{Singular unitarity in ``quantization commutes with reduction''}
\author{Hui Li}
\address{Mathematics\\
         University of Luxembourg\\
         162A, Ave de la Faiencerie\\
         L-1511, Luxembourg}
\email{li.hui@uni.lu}
$\footnote{2000 classification: 53D50, 53D20, 81S10.}$
\begin{abstract}
Let $M$ be a connected compact quantizable K\"ahler manifold equi-
pped with a Hamiltonian action of a connected compact Lie group $G$.
Let $M//G=\phi^{-1}(0)/G=M_0$ be the symplectic quotient at value
$0$ of the moment map $\phi$. The space $M_0$ may in general not be
smooth.  It is known that, as vector spaces, there is a natural
isomorphism between the quantum Hilbert space over $M_0$ and the
$G$-invariant subspace of the quantum Hilbert space over $M$. In
this paper, without any regularity assumption on the quotient $M_0$,
we discuss the relation between the inner products of these two
quantum Hilbert spaces under the above natural isomorphism; we
establish asymptotic unitarity to leading order  in Planck's
constant of a modified map of the above isomorphism under a
``metaplectic correction'' of the two quantum Hilbert spaces.
\end{abstract}
 \maketitle
 \section{Introduction}

  Let $M$ be an integral connected compact K\"ahler manifold with symplectic
  form $\omega$. Then $M$ is quantizable, i.e., there is a Hermitian
  holomorphic line bundle $L$ over $M$ with connection
  whose curvature is $-i\omega$. We consider the $k$th tensor power $L^{\otimes k}$
  of $L$.  The Hermitian structure on $L$ induces a Hermitian structure on
  $L^{\otimes k}$. The Hermitian structure on $L^{\otimes k}$  naturally
       equips the space of holomorphic sections of $L^{\otimes k}$ over $M$
       with an inner product. For each $k$, the quantum Hilbert
       space $\mathcal H(M, L^{\otimes k})$ is the space  of holomorphic sections of
    $L^{\otimes k}$ over $M$ with the inner product.

   Now, let $G$ be a connected compact Lie group acting on $M$ holomorphically
   and in a
   Hamiltonian fashion with equivariant moment map $\phi$.  Let
   $M//G=\phi^{-1}(0)/G=M_0$ be the $\mathbf{reduced\, space}$ at
   value $0$.

        Let us first consider the case when the action of $G$ on
   $\phi^{-1}(0)$ is free. Then $M_0$ is a smooth connected compact
   K\"ahler manifold. Assume that the $G$ action
   lifts to $L$ preserving the Hermitian metric. The Hermitian line bundle $L^{\otimes k}$ naturally
   descends
    to a
   Hermitian
   line bundle  $(L^{\otimes k})_0=(L^{\otimes k}|_{\phi^{-1}(0)})/G $ over $M_0$.
   The Hermitian structure on  $(L^{\otimes k})_0$  naturally
       equips
       the space of holomorphic sections of $(L^{\otimes k})_0$ over $M_0$ with an inner product.
   For each $k$, the quantum Hilbert space
    $\mathcal H(M_0,(L^{\otimes k})_0)$ is the space of holomorphic sections of
   $(L^{\otimes k})_0$ over $M_0$ with the inner product. This  is  the
   first ``reducing'' and then ``quantizing'' Hilbert space. The
    first ``quantizing'' and
  then ``reducing''  quantum Hilbert space is the $G$-invariant subspace
  $\mathcal H(M, L^{\otimes k})^G$ of $\mathcal H(M, L^{\otimes k})$.
    By Guillemin and Sternberg (\cite{GS}), there is a natural
      invertible linear map $A_k$ between $\mathcal H(M, L^{\otimes k})^G$ and
      $\mathcal H(M_0, (L^{\otimes k})_0)$. Let us call this linear map the
      Guillemin-Sternberg map. For
      quantum mechanics, the inner products of
       the quantum Hilbert spaces are also important.
  A few authors have observed that
  the  Guillemin-Sternberg map is not unitary, and, it does not become asymptotically unitary as
  $k\rightarrow\infty$.
  Moreover, they identified
  the volume of the
   $G$-orbits in the zero level set as an obstruction to asymptotic unitarity.
  We refer to the work of Flude (\cite{F}), Paoletti (\cite{Pao}),
  Ma-Zhang (\cite{MZ1} and \cite{MZ2}), Charles (\cite{Cha}) and Hall-Kirwin (\cite{HK}).
  Flude was the first who gave a formal computation of the leading-order
  term of the asymptotic density function (the function which relates the norm of an invariant holomorphic section upstairs and the norm of the descended section downstairs)
  and who obtained the non-unitarity result. Paoletti proved
  this result in his study of the asymptotic expansion of the Szeg\"o kernels upstairs and
  downstairs (using microlocal analysis (\cite{BG}, \cite{BS})). Ma and Zhang obtained
  this result (Theorem 0.10 for $E=\mathbb C$ in \cite{MZ2}) on their way
  of
  studying  the asymptotic expansion of the $G$-invariant Bergman kernel of the spin$^c$
  Dirac operator associated with vector bundles on a symplectic
  manifold.  Charles obtained this result in his study of (invariant) Toeplitz operators on $M$ and
  on the symplectic quotient $M_0$ (for torus actions) by looking at
  the relations he obtained of the (principal) symbols of the Toeplitz operators on $M$ and of
  the Toeplitz operators on
  $M_0$.
 In the recent study of the inner products of the quantum Hilbert spaces by Hall and Kirwin (\cite{HK}),
  they proved again non-unitarity by
  writing down an exact expression for the norm of an invariant holomorphic section upstairs as an integral
 over $M_0$ and by estimating the leading term of the asymptotic
 behavior of the density function.
 Moreover, for this ``free action'' case, they obtained asymptotic unitarity results for a modified
 quantization procedure. More precisely, they took the tensor products of the
        line bundles ${L^{\otimes k}}$'$s$ with the square root of the
       canonical bundle of $M$ (assuming it exists), called the metaplectic
       correction, and they  showed that a new defined
       Guillemin-Sternberg type map $B_k$ between the new quantum Hilbert spaces is invertible for all
       sufficiently large $k$, and that this map is asymptotically
       unitary to leading order as $k\rightarrow \infty$.

    In general, the action of $G$ on
   $\phi^{-1}(0)$ may not be free. Consequently, the quotient $M_0$ may  not be smooth.
   By \cite{SL} and by \cite{S}, $M_0$ is in general a stratified
   K\"ahler space, with the stratification being given by orbit types of the
   action. When there is only one orbit type, $M_0$ is still a
   smooth  K\"ahler manifold.
   In this general case when the action of $G$ on
   $\phi^{-1}(0)$ may not be free, the Hermitian line bundle $L^{\otimes k}$ descends to
   a  Hermitian V-line bundle  $(L^{\otimes k})_0$ over $M_0$. Let $\mathcal H(M_0,
    (L^{\otimes k})_0)$ still be the space of holomorphic sections
    of the V-line bundle  $(L^{\otimes k})_0$ over $M_0$ with the
    induced inner product.
   By Sjamaar (see Theorem~\ref{isom}), there is a natural linear isomorphism $A_k'$
     between $\mathcal H(M, L^{\otimes k})^G$ and $\mathcal H(M_0,
    (L^{\otimes k})_0)$.

    When the action of $G$ on
   $\phi^{-1}(0)$ is not free, the volume of the $G$-orbits in  $\phi^{-1}(0)$ is of course less
   ``uniform''. 
 One guesses by the above authors' results that  $A_k'$ would 
    not be unitary or asymptotically unitary after suitable quantum norms are defined. 
 In this paper,
   we drop the assumption that the action of $G$ on
   $\phi^{-1}(0)$ is free. We give a formula on the relation of the quantum
   norm of an invariant holomorphic section
   upstairs and the quantum norm of the descended section downstairs under the map  $A_k'$ and we give
   an asymptotic
   formula of this to leading order as $k\rightarrow\infty$. We see that  $A_k'$ is
    not unitary  and it is not asymptotically unitary.
    We still
   consider the ``metaplectic correction''.
     We give a description on how the square root of the canonical bundle of
     $M$ descends to $M_0$, we show the existence of a family of modified
     isomorphisms $B_k'$ for sufficiently large  $k$ between the new quantum Hilbert spaces, and we
     establish   asymptotic unitarity to leading order term for
     the maps $B_k'$.\\

     There are two main problems that need to be addressed in this
     new study. One is that we find a suitable way to descend the
     half form bundle of $M$ to the stratified quotient $M_0$. Another
     problem has to do with a large piece of the manifold $M$, the
     semistable set $M^{ss}$, which is open dense and connected in
     $M$.  By Theorem~\ref{isom}, $\mathcal H(M_0, (L^{\otimes k})_0)\simeq\mathcal H(M^{ss}, L^{\otimes
     k})^G$.   The holomorphic action of $G$ can be analytically
     extended to a $G_{\mathbb C}$-action, where $G_{\mathbb C}$ is the
     complexification of $G$.
     If $G$ acts freely on $\phi^{-1}(0)$,
     $M^{ss}$  consists of free $G$-orbits and it consists of
     complex $G_{\mathbb C}$-orbits each of which  intersects
     $\phi^{-1}(0)$ at one $G$-orbit. In the general case,
     $M^{ss}$ may contain
     complex $G_{\mathbb C}$-orbits which do not intersect
     $\phi^{-1}(0)$ but contain those $G_{\mathbb C}$-orbits which intersect $\phi^{-1}(0)$
     in their closures. We will analyze the structure of these complex orbits and study their contribution
     to the quantum norms. \\

     Our main results are Theorem~\ref{redhfb'}, Theorem~\ref{modisom}, Theorem~\ref{qnsr},
     Theorem~\ref{asym} (and the corollaries of Theorems~\ref{qnsr}
     and \ref{asym},
     corollaries~\ref{qnsr1} and \ref{qnsr2}),
     and Theorem~\ref{uni}.\\

     We will use three different notations interchangeably for the
     symplectic quotient at $0$, $M_0$, $M//G$, and
     $M^{ss}//G_{\mathbb C}$, depending on the context.

\subsubsection*{Acknowledgement}

     This work was mainly motivated by the work \cite{HK} of Brian Hall and William
  Kirwin. I thank them for their work.
   I thank Reyer Sjamaar for answering me an email inquiry on semistable points when
  I started to consider the problem. I thank Laurent Charles, Xiaonan Ma and Weiping
  Zhang, and Roberto Paoletti for pointing to me their articles.
 I was happy
     to see through their work
 related problems and areas. Finally, I am grateful to the referee
 for making comments and remarks which helped me to improve the
 exposition.

  \section{Reduction of K\"ahler manifolds}

   In this section, we will recall some main results
   obtained by R. Sjamaar in \cite{S} on general K\"ahler quotients.
   This will help us to understand our space $M_0$ as well as its relation with $M$.
   One may see  the difference between the case when the action of $G$ on
   $\phi^{-1}(0)$ is free and the case when this assumption is
   removed.
   This section  also serves as a
   preparation for the tools needed in the subsequent sections.

    To understand stratified K\"ahler spaces and their quantizations, we also refer
    to the work of Huebschmann,
    \cite{H} and \cite{H'}.\\

    Let $(M, \omega, J, B=\omega(\cdot, J\cdot))$ be a connected compact K\"ahler manifold with symplectic
    form $\omega$, compatible complex structure $J$ and Riemannian metric $B$. Let $G$ be a connected compact
    Lie group acting holomorphically on $M$. Assume that the $G$ action is Hamiltonian with an equivariant moment
    map $\phi$. Assume $a$ is a value of $\phi$. Then the quotient $M_a=\phi^{-1}(G\cdot a)/G$ is called
    the symplectic quotient or the reduced space at the coadjoint orbit $G\cdot a$. Let us restrict attention
    to the value $a=0$. By \cite{SL}, the quotient $M_0$ is a connected compact
   stratified symplectic space with a connected open dense  stratum. If $M_0$ has only one stratum, then it is a smooth
   symplectic manifold. We will see that $M_0$ also admits an analytic
   structure such that $M_0$ is a stratified K\"ahler space.

    Since $G$ acts holomorphically, the action can be analytically
    continued to a holomorphic action of the ``complexified'' group $G_{\mathbb{C}}$ on $M$. The Lie algebra
    $\mathfrak g_{\mathbb{C}}$
    of $G_{\mathbb{C}}$ is the complexification of $\mathfrak g$. The Cartan decomposition gives a diffeomorphism
    $G_{\mathbb{C}}\backsimeq exp(i\mathfrak g)G$. For $\xi\in \mathfrak g$, let $X^{\xi}$ be the infinitesimal vector
    field on $M$ generated by $\xi$. Then $X^{i\xi}=JX^{\xi}$ is the infinitesimal vector field generated
    by $i\xi$.

     Define a point $m$ in $M$ to be $\mathbf{(analytically)\, semistable}$ if the closure of the  $G_{\mathbb{C}}$-orbit through $m$ intersects the zero level set $\phi^{-1}(0)$. Let $M^{ss}$ be the set of semistable points
     in $M$.   The point $m$ is called
    $\mathbf{(analytically)\, stable}$ if the closure of the  $G_{\mathbb{C}}$-orbit through $m$ intersects the zero level set $\phi^{-1}(0)$ at a point where $d\phi$ is surjective. Let
      $M^{s}$ be the set of stable points in $M$. When the action of
      $G$ on $\phi^{-1}(0)$ is free or locally free, $d\phi$ is surjective at any
      point of $\phi^{-1}(0)$. In this case,  $M^{ss}$ coincides
      with $M^{s}$.

    Assuming there is a $G$-invariant inner product on $\mathfrak
    g$, by Lemma 6.6 in \cite{K}, the gradient vector field of $\|\phi\|^2$
    is given by
    $$grad(\|\phi\|^2)(m)=2JX^{\phi(m)}(m),$$
    where we have identified $\phi(m)\in\mathfrak g^*$ with a vector
    in $\mathfrak g$ using the inner product, and where $X^{\phi(m)}(m)$ is the vector field on $M$
    induced by $\phi(m)$, evaluated
    at the point $m$. So $grad(\|\phi\|^2)(m)$ is tangent to the
    $G_{\mathbb{C}}$-orbits. Let $F_t$ be the flow of
    $-grad(\|\phi\|^2)$. Kirwan has proved that $M^{ss}$ is the set of
    points $m\in M$ such that the path $F_t(m)$ has a limit point in
    $\phi^{-1}(0)$ (\cite{K}). By \cite{L} or by \cite{W}, the limit
    map $F_{\infty}(m)$ gives an equivariant deformation retraction
    from $M^{ss}$ onto  $\phi^{-1}(0)$.
 \subsection{The holomorphic slice theorem}
  In order to describe the complex analytic structure on $M_0$, let us first recall
  the holomorphic slice theorem due to R. Sjamaar. The results on the orbit structure of $M^{ss}$ and on
  the stratified  K\"ahler structure of $M_0$ are due to this theorem.

  \begin{theorem}\label{slice}
  (Holomorphic slice theorem) (\cite{S})
  Let $M$ be a K\"ahler manifold and let  $G_{\mathbb{C}}$ act holomorphically
  on $M$. Assume the action of the compact real form $G$ is Hamiltonian. Let $m$ be any
  point in $M$ such that the $G$-orbit through $m$ is isotropic. Then there exists a $\mathbf{holomorphic\, slice}$
 at $m$ for the $G_{\mathbb{C}}$-action.
\end{theorem}

If $X$ is a complex space and  $G_{\mathbb{C}}$  a reductive complex
Lie group
 acting holomorphically on $X$, we have the following definition of a
 $\mathbf{holomorphic\, slice}$.

  \begin{definition}
 A $\mathbf{holomorphic\, slice}$ at $x$ for the $G_{\mathbb{C}}$ action is a locally closed analytic subspace
 $D$ of $X$ with the following properties:\\
 1. $x\in D$;\\
 2. $G_{\mathbb{C}}D$ of $D$ is open in $X$;\\
 3. $D$ is invariant under the action of the stabilizer $(G_{\mathbb{C}})_x$;\\
 4. the natural $G_{\mathbb{C}}$-equivariant map from the associated bundle
     $G_{\mathbb{C}}\times_{(G_{\mathbb{C}})_x} D$ into $X$, which sends an
   equivalence class $[g, y]$ to the point $gy$, is an analytic isomorphism onto
  $G_{\mathbb{C}}D$.
  \end{definition}

  \subsection{K\"ahler reduction}

     For the orbit structure of $M^{ss}$, one may see Proposition 2.4 in
     \cite{S}. We list a few of them which are more relevant to us.
     \begin{proposition}\label{propS}
      In the following, ``closed'' means ``closed in $M^{ss}$'' and
      ``closure'' means ``closure in $M^{ss}$''.

      1. The semistable set $M^{ss}$ is the smallest   $G_{\mathbb{C}}$-invariant open subset of $M$ containing
      $\phi^{-1}(0)$, and its complement in $M$ is a complex-analytic
      subset;

      2. A $G_{\mathbb{C}}$-orbit in $M^{ss}$ is closed if and only
      if it intersects $\phi^{-1}(0)$;

      3. The closure of every   $G_{\mathbb{C}}$-orbit in $M^{ss}$
      contains exactly one closed  $G_{\mathbb{C}}$-orbit.

     \end{proposition}

     We call two semistable points $x$ and $y$ $\mathbf{related}$ if the closures in  $M^{ss}$
      of the orbits  $G_{\mathbb{C}}x$
      and  $G_{\mathbb{C}}y$ intersect. This relation is an equivalence relation. Let $M^{ss}//G_{\mathbb{C}}$
     be the quotient space and let $\pi_{\mathbb{C}}: M^{ss}\rightarrow M^{ss}//G_{\mathbb{C}}$ be the quotient map.

      \begin{theorem} (\cite{S})
      The inclusion $\phi^{-1}(0)\subset M^{ss}$ induces a homeomorphism
      $M_0=\phi^{-1}(0)/G\rightarrow  M^{ss}//G_{\mathbb{C}}$.
      \end{theorem}

     We say that a subset $A$ of $M^{ss}$ is $\mathbf{saturated}$
     with respect to $\pi_{\mathbb C}$ if $\pi_{\mathbb C}^{-1}\pi_{\mathbb
     C}(A)=A$.

    \begin{proposition}\label{saturated}(\cite{S})
     At every point of $\phi^{-1}(0)$, there exists a
     $\mathbf{holomorphic\, slice}$ $D$ such that the set $G_{\mathbb
     C}D$ is $\mathbf{saturated}$ with respect to the quotient mapping $\pi_{\mathbb
     C}$.
    \end{proposition}

      We identify the spaces $M^{ss}//G_{\mathbb{C}}$ and $M_0$. We furnish $M_0$ with a complex-analytic structure such that the
      quotient map $\pi_{\mathbb{C}}$ is holomorphic. We define a function $f$ defined on an open subset $O$ of $M_0$ to be holomorphic if
      the pullback of $f$ to $\pi_{\mathbb{C}}^{-1}(O)$ is holomorphic. Let $\mathcal O_{\mathcal M_0}$ be the sheaf of
      holomorphic functions on $M_0$.

       \begin{theorem}\label{analytic} (\cite{S})
       The ringed space $(M_0,\mathcal O_{\mathcal M_0})$ is an analytic space.
        \end{theorem}

     The following theorem describes the property of the stable set
     $M^s\subset M^{ss}$. If $0$ is a regular value of $\phi$, then $M^s=M^{ss}$.
     In general, if  $M^s\neq\emptyset$, then $M^s$ is  open and dense in $M^{ss}$.
      \begin{theorem}(\cite{S})\label{theorem-stable}
      If $x\in M$ is stable, then the orbit  $G_{\mathbb{C}}x$ is closed
      in $M^{ss}$. Let $Z$ be the set of points  $m\in\phi^{-1}(0)$
      with the property that $d\phi_m$ is surjective. Then the
      stable set $M^s$ is equal to $F_{\infty}^{-1}(Z)$. Every fiber
      of $\pi_{\mathbb{C}}|_{M^s}$ consists of a single orbit.
      \end{theorem}

      By this theorem, we see that if a $G$-orbit $\mathrm O=G\cdot x$ in $\phi^{-1}(0)$
      has the dimension of $G$, then only one complex orbit
      $G_{\mathbb{C}}\cdot x=G_{\mathbb{C}}\cdot\mathrm O$ flows to
      $\mathrm O$ under the gradient flow of $-\|\phi\|^2$.\\

      The stratification of
      $M_0$ as a stratified
      symplectic space  is given by orbit types.
      Let $p\in M_0$, and let $x\in\pi^{-1}(p)$, where $\pi: \phi^{-1}(0)\rightarrow M_0$
      is the quotient map. Let $(H)$ be a conjugacy
      class of closed subgroups of $G$. Then $p$ is said to be of
      orbit type $(H)$ if the stabilizer of $x$ is conjugate to $H$.
      By \cite{SL}, the set of all points of orbit type $(H)$ in $M_0$ is a
      symplectic manifold.

      We can similarly define $G_{\mathbb{C}}$-orbit types. We can
      show that if $x\in\phi^{-1}(0)$, then the complex stabilizer
      $(G_{\mathbb{C}})_x$ is equal to the complexification
      $(G_x)_{\mathbb{C}}$ of the compact stabilizer $G_x$ (see Proposition 1.6 in \cite{S}).
       By
      Proposition~\ref{propS}, the fiber $\pi_{\mathbb{C}}^{-1}(p)$
      contains a unique closed  $G_{\mathbb{C}}$-orbit
      $G_{\mathbb{C}}x$. Let us say $p$ is of $G_{\mathbb{C}}$-orbit
      type $(H_{\mathbb{C}})$ if the stabilizer $(G_{\mathbb{C}})_x$
      is conjugate to $H_{\mathbb{C}}$ in $G_{\mathbb{C}}$.

      \begin{theorem}\label{strkahler} (\cite{S})
      The stratification of $M_0$ by $G$-orbit types is identical to
      the stratification by $G_{\mathbb{C}}$-orbit types. Each
      stratum $\mathcal S$ is a complex manifold and its closure is a
      complex-analytic subvariety of $M_0$. The reduced symplectic
      form on $\mathcal S$ is a K\"ahler form.
      \end{theorem}

  \section{quantization of K\"ahler manifolds}
   Let $M$ a  connected compact K\"ahler manifold as in the last
   section. Assume that the K\"ahler form $\omega$ is integral, i.e.,
   the cohomology class $[\omega/2\pi]$ is an integral cohomology
   class. Then $M$ is quantizable, i.e., there is a Hermitian line
   bundle $L$ with compatible connection $\triangledown$ such that
   its curvature is $-i\omega$. The $k$-th tensor power $L^{\otimes k}$ of $L$
   is a Hermitian line bundle over $M$ with induced Hermitian
   structure from $L$. For each $k$, $L^{\otimes k}$ may be given
   the structure of a holomorphic line bundle. For each fixed
   $k$, the quantum Hilbert space is the space of holomorphic
   sections of $L^{\otimes k}$ over $M$, denoted $\mathcal H(M, L^{\otimes
   k})$. Let $\epsilon_{\omega}=\frac{\omega^n}{n!}$ be the Liouville volume
   form on $M$. Then the inner product on $\mathcal H(M, L^{\otimes k})$ is usually defined
   to be
        $$<s_1, s_2>=(k/2\pi)^{n/2}\int_M (s_1, s_2)\epsilon_{\omega},$$
   where $(s_1, s_2)$ is the pointwise Hermitian structure on
   $L^{\otimes k}$.

       In this paper, we will study quantizable K\"ahler
       manifolds with a holomorphic Hamiltonian Lie group action.
       The symplectic quotient at value $0$ may not be smooth.
      To adapt to this situation, we will give two
       definitions of the inner product on $\mathcal H(M, L^{\otimes
       k})$, respectively in Definitions~\ref{def-qn1} and \ref{def-qn2} of Section 13.

  \section{quantum reduction}
   Let $M$ be a connected compact quantizable K\"ahler
   manifold. Let $G$ be a connected compact Lie
    group acting on $M$
    holomorphically and in a Hamiltonian fashion with moment map $\phi$.
   The $G$ action lifts to a
    holomorphic action on the line bundle $L$ preserving the
    Hermitian structure. Both the $G$ action on $M$ and on $L$ can
    be analytically continued to holomorphic $G_{\mathbb{C}}$
    actions. The $G$ action on $L$ induces $G$ actions on $\mathcal H(M, L^{\otimes
    k})$. Infinitesimally, the action is given by
    $$Q_{\xi}s=\triangledown^{(k)}_{X^{\xi}}s-ik\phi_{\xi}s,\,\, \mbox{for}\,\,\xi\in\mathfrak g,$$
    where $\triangledown^{(k)}$ is the induced connection on  $L^{\otimes k}$, and $\phi_{\xi}$ is the ``$\xi$-moment map component", i.e., $\phi_{\xi}=<\phi, \xi>$.
    The reduction at quantum level amounts to taking $G$-invariant
    holomorphic sections, i.e., taking $\mathcal H(M, L^{\otimes k})^G$.

   \section{quantization after reduction}
      Let $M$ be a connected compact quantizable K\"ahler manifold
   equipped with a holomorphic Hamiltonian action of a connected compact Lie
    group $G$. Let $\phi$ be the moment map. Let $L_0=L|_{\phi^{-1}(0)}/G$.
   Then $L_0$ is a $V$-line bundle over $M_0$, i.e., each point in $M_0$ has an open neighborhood
   $O$ which is the quotient of a space $\tilde O$ by a finite group $\Gamma$ such that $L_0|_O$ is
   the quotient by $\Gamma$ of a $\Gamma$-equivariant line bundle over $\tilde O$.  As an analytic space, $L_0$
   can be identified with the quotient $L|_{M^{ss}}//G_{\mathbb{C}}$. A holomorphic section
   of $L$ defined over a $G_{\mathbb{C}}$-invariant open set is $G$-invariant if
   and only if it is $G_{\mathbb{C}}$-invariant. Let $\mathcal L$ be the sheaf of holomorphic sections of $L$ and
    define a sheaf $\mathcal L'$ on $M_0$, the sheaf of invariant sections, by letting
    $\mathcal L'(O)=\mathcal L(\pi_{\mathbb{C}}^{-1}(O))^{G_{\mathbb{C}}}$ for each open set $O$ of $M_0$.
    Then we have
    \begin{proposition} (\cite{S})\label{propS'}
     The sheaf $\mathcal L'$ is (the sheaf of sections of) the holomorphic $V$-line bundle $L_0$
     over $M_0=M^{ss}//G_{\mathbb{C}}$.
    \end{proposition}

    We take the space of holomorphic sections $\mathcal H(M_0, L_0)$ of the $V$-line bundle $L_0$ as the
    quantization of the reduced space $M_0$.

     If we replace $L$ by $L^{\otimes k}$, we have the quantum spaces $\mathcal H(M_0, (L^{\otimes k})_0)$.

     Since the action of $G$ preserves the Hermitian structure on $L^{\otimes k}$, the Hermitian structure
     on $L^{\otimes k}$ descends to a Hermitian structure on $(L^{\otimes k})_0$. Let $s\in\mathcal H(M, L^{\otimes k})^G$. Then by restricting $s$ to
     $\phi^{-1}(0)$ and by letting it descend to $M_0$, we get an element of $\mathcal H(M_0, (L^{\otimes k})_0)$. Let
     us call this linear map $A_k'$. So, if $x\in\phi^{-1}(0)$, then we have
      $$|s|^2(x)=|A_k' s|^2([x]).$$

     We still need to define an inner product on  $\mathcal H(M_0, (L^{\otimes k})_0)$.
   Denote
     \begin{equation}\label{Z_H}
 Z_{(H)}=\{m\in\phi^{-1}(0): \mbox{the stabilizer group of}\,\, m \,\,\mbox{is conjugate to}\,\, H\subset G
 \},
    \end{equation}
    and
     \begin{equation}\label{S_H}
     \mathcal S_{(H)}=Z_{(H)}/G.
     \end{equation}

      As Remarked by Sjamaar (see Remark 3.9 in \cite{SL}),
    each  $\mathcal S_{(H)}$ has finite symplectic volume. For a fixed  stratum
     $\mathcal S$ of
     $M_0$, let $d_{\mathcal S}$ be the complex dimension of
   $\mathcal S$, and let $\epsilon_{\hat{\omega}_{\mathcal S}}$ be the volume form on $\mathcal S$,
     where $\hat{\omega}_{\mathcal S}$ is the reduced symplectic
     form on $\mathcal S$.

     Let $s_1', s_2'\in \mathcal H(M_0, (L^{\otimes k})_0)$, and,
     let $(s_1',
     s_2')$ be the pointwise Hermitian inner product on $(L^{\otimes k})_0$
     inherited from the one on $L^{\otimes k}$.
     Since there is an open dense connected stratum, say
         $\mathcal S^O$, in $M_0$, which has full measure, we give the first definition of
         an inner product on  $\mathcal H(M_0, (L^{\otimes k})_0)$:

         \begin{equation}\label{norm_1}
     <s_1', s_2'>_{(1)}=(k/2\pi)^{d_{\mathcal S^O}/2}\int_{\mathcal S^O} (s_1', s_2')\epsilon_{\hat{\omega}_{\mathcal S^O}}.
     \end{equation}

      The following second definition  of an inner product on  $\mathcal H(M_0, (L^{\otimes k})_0)$
      takes into account all the strata of
      $M_0$:
     \begin{equation}\label{norm_2}
     <s_1', s_2'>_{(2)}=\sum_{\mathcal S_{(H)}}(k/2\pi)^{d_{\mathcal S_{(H)}}/2}\int_{\mathcal S_{(H)}} (s_1', s_2')\epsilon_{\hat{\omega}_{\mathcal
     S_{(H)}}}.
     \end{equation}

     If $\mathcal S$ is a
     single point, then the above integral of $(s_1',
     s_2')$  over $\mathcal S$ is just the value of $(s_1',
     s_2')$ over this point.

  \section{The linear space isomorphism}
     In the last section, we defined a linear map $A_k'$ from $\mathcal H(M, L^{\otimes
     k})^G$ to $\mathcal H(M_0, (L^{\otimes k})_0)$. We have

       \begin{theorem}\label{isom} (\cite{S})
    Under our hypotheses, the quotient map $\pi_{\mathbb{C}}: M^{ss}\rightarrow
    M_0$ and the inclusion $M^{ss}\subset M$ induce isomorphisms
    $\mathcal H(M_0, (L^{\otimes k})_0)\simeq\mathcal H(M^{ss}, L^{\otimes
     k})^G\simeq\mathcal H(M, L^{\otimes k})^G$.
    \end{theorem}

    By Proposition~\ref{propS'}, we have the isomorphism   $\mathcal H(M_0, (L^{\otimes k})_0)\simeq\mathcal H(M^{ss}, L^{\otimes
     k})^G$. The isomorphism $\mathcal H(M^{ss}, L^{\otimes
     k})^G\simeq\mathcal H(M, L^{\otimes k})^G$ is based on the observation that the norm
    of an invariant holomorphic section $s$ of $L^{\otimes k}$ is
    increasing along the trajectories of $-grad(\|\phi\|^2)$. It
    follows that if $s$ is defined on $M^{ss}$, then $<s, s>$ is
    bounded on $M$. By Riemann's Extension Theorem, $s$ extends to
    a $G$-invariant holomorphic section on $M$.  See \cite{S} for details.\\

    From this theorem, we can deduce that a point $x\in M$ is semistable if
 there exists
 an invariant global holomorphic section $s\in\mathcal H(M, L^{\otimes l})^G $ for
 some  $l$ such that $s(x)\neq 0$ (see \cite{S}). So the set of unsemistable points
 is contained in the $0$ set of $s$, therefore it has complex codimension at least one.

     \section{Half form bundles on $M$}
     Let $K=\bigwedge^n(T^{1, 0}M)^*$ be the canonical
     bundle of $M$. A smooth section of $K$ is called an $(n, 0)$-form.
     We know that the first Chern class of $K$ is $-c_1(M)$.
        Assume $c_1(M)/2$ is integral. Then the square root $\sqrt K$ of the bundle $K$ exists.
        We fix a choice of  $\sqrt K$. The group $G$ acts on
        sections of $K$. 
        Infinitesimally,
        a Lie algebra element $\xi\in \mathfrak g$ acts on $(n,
        0)$-forms by taking the Lie derivative $L_{X^{\xi}}$ of the form.
        This induces an action of $\mathfrak g$ on half forms by
       $2(L_{X^{\xi}}\mu)\mu=L_{X^{\xi}}(\mu^2)$,
       where $\mu\in\sqrt K$. Since $G$ acts holomorphically
       on $M$, we can check that $\mathfrak g$ preserves the space of
       holomorphic sections of  $K$ and of $\sqrt K$.

            Let us define a Hermitian structure on
       $\Gamma(M, \sqrt{K})$, where  $\Gamma(M, \sqrt{K})$ is the space of smooth sections
       of $\sqrt{K}$. Let $\mu, \nu\in \Gamma(M, \sqrt{K})$ be half
       forms, then
       $\mu^2\wedge\bar{\nu}^2\in\Gamma(\bigwedge^{2n}T^*(M))$. The
       volume form $\epsilon_{\omega}$ is a global trivializing section of
       $\bigwedge^{2n}T^*(M)$. So there is a function, denoted
       $(\mu, \nu)$, such that
     \begin{equation}\label{hermitian-srb}
       \mu^2\wedge\bar{\nu}^2=(\mu, \nu)^2\epsilon_{\omega}.
     \end{equation}
       The function  $(\mu, \nu)$ is defined to be the pointwise inner
       product of $\mu$ and $\nu$.

      We use this to define a Hermitian form on $\Gamma(M, L^{\otimes k}\otimes\sqrt{K})$.
   Let $t_1, t_2\in\Gamma(M, L^{\otimes k}\otimes\sqrt{K})$ which are locally represented
   by $t_j(x)=s_j(x)\mu_j(x)$, we define
   \begin{equation}\label{hermitian-product}
   (t_1, t_2)(x)=(s_1(x), s_2(x))(\mu_1, \mu_2)(x).
   \end{equation}

   \section{The ``push down" of the half form bundle $\sqrt{K}$ to the reduced space $M_0$}

  \subsection{When the $G$ action on $\phi^{-1}(0)$ is free}

   Before we come to the general case, let us recall first the
   procedure given by Hall and Kirwin of pushing down  a half form bundle $\sqrt K$ of $M$
   to a half form bundle $\sqrt{\hat K}$ on $M//G$ in the case the
   action of $G$ on $\phi^{-1}(0)$ is free.\\

  Let  $\alpha$ be a $G_{\mathbb{C}}$-invariant $(n, 0)$-form on $M$.
  Hall and Kirwin obtained an
   $(n-d, 0)$-form  ($d$ is the dimension of $G$) $\hat\beta$ on $M//G$
   in the following way. Choose a $G$-invariant inner product on $\mathfrak g$
normalized so
 that the volume of $G$ with respect to the associated Haar measure is 1. Fix an orthonormal basis ${ \xi_1, \xi_2, ... , \xi_d}$ of the
  Lie algebra $\mathfrak g$. Let ${X^{\xi_1}, X^{\xi_2}, ... , X^{\xi_d}}$ be
  the vector fields they generate on $M$. For any $x\in M^s$, define

  $$\beta=i(\bigwedge_jX^{\xi_j})\alpha.$$

  One can show that $\beta$ is  basic with respect to the projection map
  $\pi_{\mathbb{C}}$. So  $\beta=\pi_{\mathbb{C}}^*(\hat\beta)$, where
  $\hat\beta$ is an $(n-d, 0)$-form on $M//G$. Let $\mathfrak B$
  be the map
     $$\mathfrak B(\alpha)=\hat\beta.$$
     Conversely, one can construct the inverse map of this push down map. Given an  $(n-d,
 0)$-form  $\hat\beta$ on $M//G$. The pull back $\beta=\pi_{\mathbb{C}}^*(\hat\beta)$
 is a  $G_{\mathbb{C}}$-invariant  $(n-d,
 0)$-form on $M^s$. One can construct a $G_{\mathbb{C}}$-invariant $(n, 0)$-form $\alpha$ on $M^s$ from
 $\beta$. Given a local frame ${X^{\xi_1}, X^{\xi_2}, ..., X^{\xi_d}, Y_1,
 ..., Y_{n-d}}$ for $T_xM^s$, set
  $$\alpha (X^{\xi_1}, X^{\xi_2}, ..., X^{\xi_d}, Y_1,..., Y_{n-d})=\pi_{\mathbb{C}}^*\hat\beta(Y_1, ...,
   Y_{n-d}),$$
  and define $\alpha$ on any other frame by $GL(n,
\mathbb{C})$-equivariance
 and the requirement that $\alpha$ be an $(n, 0)$-form. Every other frame is equivalent to
a linear combination of frames which are $GL(n,
\mathbb{C})$-equivalent to one of the form ${W_1, W_2, ..., W_d, Y_1,
 ..., Y_{n-d}}$ where $W_j=X^{\xi_j}$ or $JX^{\xi_j}$.

   Assume that the $\mathfrak g$ action on $\sqrt K$ exponentiates to a $G$ action and
 it is compatible with the $G$ action on $K$. It can be shown that the $G$ action on
 $\sqrt K$ can be analytically continued to a  $G_{\mathbb{C}}$ action. Define a line
 bundle $\sqrt{\hat K}$ over $M//G$ whose fiber is the equivalence class of $\sqrt K$ under the
 $G_{\mathbb{C}}$ action.
  For a $G_{\mathbb{C}}$-invariant smooth section $\mu\in\Gamma(M, \sqrt K)^{G_{\mathbb{C}}}$, we define the map
$$B: \Gamma(M, \sqrt K)^{G_{\mathbb{C}}}\rightarrow \Gamma(M//G, \sqrt{\hat
   K})$$
  by $(B\mu)^2=\mathfrak B(\mu^2)$.

   Since for an $(n, 0)$-form $\alpha$, contracting with $\bigwedge_j X^{\xi_j}$ is the same as
   contracting with $\bigwedge_j \pi_+X^{\xi_j}$, where $\pi_+X^{\xi_j}=\frac{1}{2}(X^{\xi_j}-iJX^{\xi_j})$,
    and the vector
   fields $\pi_+X^{\xi}$ are holomorphic, $\alpha$ is locally
   holomorphic if and only if $\mathfrak B (\alpha)$ is locally
   holomorphic; and, $\mu$ is locally holomorphic if and only if
   $B(\mu)$ is locally holomorphic.\\

\subsection{When the $G$ action on $\phi^{-1}(0)$ is not necessarily free}

  Now, we come to the general case.

 \begin{lemma}\label{pushdown}
 Let $\alpha\in\Gamma(M, K)^{G_{\mathbb{C}}}$. Then,
 $\alpha$ descends to a smooth $(d_{\mathcal S}, 0)$-form $\hat
\beta|_{\mathcal S}$ on each smooth stratum  $\mathcal S$ of $M_0$
of complex dimension $d_{\mathcal S}$. If $\alpha$ is holomorphic,
then each $\hat \beta|_{\mathcal S}$ is holomorphic.
 \end{lemma}
 \begin{proof}
 Let
 $Z_{(H)}$ and $\mathcal S_{(H)}$ be as in (\ref{Z_H}) and (\ref{S_H}).
 Take the complex submanifold $G_{\mathbb{C}}\cdot Z_{(H)}$.
 Let
  $$\alpha_|=\alpha|_{G_{\mathbb{C}}\cdot Z_{(H)}}.$$
Then $\alpha_|$ is a $G_{\mathbb{C}}$-invariant  $(m, 0)$-form on
$G_{\mathbb{C}}\cdot Z_{(H)}$,
  assuming  $m$ is the complex dimension of $G_{\mathbb{C}}\cdot
  Z_{(H)}$.\\

Case 1. Assume  $H=G$. Then
  $G_{\mathbb{C}}\cdot Z_G=Z_G=Z_G/G=\mathcal S_G$.
  We define  $\hat\beta|_{\mathcal S_G}=\alpha|_{Z_G}$.\\

Case 2. Assume  $H\neq G$.
 Assume we have chosen a normalized
$G$-invariant inner
 product on $\mathfrak g$.
Let ${\xi_1, ..., \xi_h}$ be  an orthonormal basis of $\mathfrak h=Lie(H)$, expand it to an orthonormal basis of
 $\mathfrak g=Lie(G)$ by joining ${\xi_{h+1}, ..., \xi_d}$.  At each point $x$ of $Z_{(H)}$
 with stabilizer group $H$,
 we define
  $$\beta_|=i(\bigwedge_{j=h+1, ..., d}X^{\xi^j})\alpha_|.$$
 We contract the form $\alpha_|$ similarly at the points of
 $Z_{(H)}$ with stabilizers conjugate to $H$. So, along $Z_{(H)}$, we have a new form $\beta_|$.
  Let
  $${\beta_|}_|=\beta_||_{Z_{(H)}}.$$
 This restriction ``cuts off'' the $JX^{Ad(G)\cdot\xi^j}$, $j=h+1, ..., d$
 directions  which are normal to $Z_{(H)}$ in $G_{\mathbb{C}}\cdot Z_{(H)}$.
 Now, ${\beta_|}_|$
  is a  smooth $G$-invariant $(m-d_{G/H})$-form
  defined on $Z_{(H)}$. By the above contraction and by
  $G$-invariance of the form  $\alpha$,
 clearly, ${\beta_|}_|=\pi^*(\hat{\beta_|})$, where $\hat{\beta_|}$
 is a
 $(m-d_{G/H}, 0)$-form on $\mathcal S_{(H)}$,
  $\pi: Z_{(H)}\rightarrow  \mathcal S_{(H)}$ is the quotient map.

  By the above construction, if $\alpha$ is holomorphic,
  then  $\hat \beta|_{\mathcal S}$ is holomorphic (see the reason we mentioned in Section 8.1).
 \end{proof}

  Next,  we will use the holomorphic slice theorem to see how the forms $\hat{\beta}|_{\mathcal S}$'$s$
  are related.

  \begin{lemma}\label{pushdown'}
  The forms $\hat\beta|_{\mathcal S}$'$s$ in Lemma~\ref{pushdown}
  satisfy:  if $\mathcal S\subset \bar{\mathcal S'}$, then $\hat\beta|_{\mathcal S}$
  is obtained from $\hat\beta|_{\mathcal S'}$ by degenerating some
  directions.
  \end{lemma}
 \begin{proof}
 Let $x_0\in Z_{(H)}$ be a point with stabilizer group $H$.
  Take a
  $\mathbf{saturated}$ open
  neighborhood $U=G_{\mathbb{C}}D=G_{\mathbb{C}}\times_{H_{\mathbb
  C}}D$ of $x_0$ (see Proposition~\ref{saturated}), where $H_{\mathbb C}=(H)_{\mathbb C}$ is the complex stabilizer group of
  $x_0$ which is the complexification of $H$.
  Split $D=D_1\times D_2$, where $D_1$ is the
  fixed complex subspace of the $H$ action and therefore the $H_{\mathbb C}$ action. So
  $U=G_{\mathbb{C}}\times_{H_{\mathbb
  C}}(D_1\times D_2)$. The set
    $$U//G_{\mathbb{C}}=(D_1\times D_2)//H_{\mathbb{C}}=D_1\times D_2//H_{\mathbb{C}}$$
    is a neighborhood of $[x_0]$ in $M_0$.

  The set $U$ is $G$-equivariantly diffeomorphic
 to $G\times_{H}(\sqrt{-1}\mathfrak m\times D_1\times D_2)$, where $\mathfrak m$
 is the orthogonal complement of $\mathfrak h$ in $\mathfrak g$. We
 pull back (or restrict) the symplectic form $\omega$ on $M$ to $U$.
 The group $H_{\mathbb C}$ acts on $D_2$ holomorphically.
  Assume $\phi_|$ is the moment map for the $H$-action on $D_2$ with respect to
  the restricted K\"ahler form.
  Then,
  $$U\cap\phi^{-1}(0)=G\times_H (D_1\times \phi_|^{-1}(0)).$$

  Denote
   $$Z'_{(H')}=\{m\in\phi_|^{-1}(0)\subset D_2: \mbox{the compact stabilizer group of}\,\, m\,\,\mbox{is
  conjugate to}\,\, H'\subset H\},$$
  and recall  (\ref{Z_H}) for  $Z_{(H)}$. We have
   $$U\cap Z_{(H')}=G\times_H (D_1\times Z'_{(H')}),\,\,\mbox{in particular},\,\, U\cap Z_{(H)}=G\times_H (D_1\times 0),$$
  where $0$ is in the closure of  $Z'_{(H')}$.
 While
  $$G_{\mathbb{C}}\cdot (U\cap Z_{(H')})=G_{\mathbb{C}}\times_{H_{\mathbb C}} (D_1\times H_{\mathbb C} Z'_{(H')}),
  \,\,\mbox{and},\,\,
 G_{\mathbb{C}}\cdot (U\cap Z_{(H)})=G_{\mathbb{C}}\times_{H_{\mathbb C}}(D_1\times 0).$$

  The quotients are
  $$(U\cap Z_{(H')})/G=D_1\times Z'_{(H')}/H, \,\,\mbox{which is the same as}$$
  $$G_{\mathbb{C}}\cdot (U\cap Z_{(H')})//G_{\mathbb{C}}=D_1\times (H_{\mathbb C} Z'_{(H')})//H_{\mathbb C},$$
  and
   $$(U\cap Z_{(H)})/G=D_1\times 0,\,\,\mbox{which is the same as}\,\,  G_{\mathbb{C}}\cdot (U\cap Z_{(H)})//G_{\mathbb{C}}.$$

  Now, restricting to the open set $U$, resp., $U//G_{\mathbb{C}}$, the
  relation between $Z_{(H)}$ and $Z_{(H')}$, resp., the relation
  between $\mathcal S_{(H)}$ and $\mathcal S_{(H')}$ is clear. In the open set $U$, do
  the specified restricting, contracting, restricting again, and pushing down of the form
 $\alpha_|$ as in the proof of Lemma~\ref{pushdown}, we see that, if
 we use local coordinates, and if $\hat\beta|_{\mathcal S_{(H')}}=g(w_1, ..., w_k) dw_1\wedge...\wedge dw_k$,
 then $\hat\beta|_{\mathcal S_{(H)}}=g(w_{i_1}, ..., w_{i_j}, 0,..., 0)dw_{i_1}\wedge...\wedge dw_{i_j}$,
 where $\{i_1,..., i_j\}\subset\{1,..., k\}$.
 \end{proof}

  Let us simply use
  $\hat\beta$ to denote this family
 of forms on $M_0$ we got. It has different dimensions on different dimensional strata.

   Let us denote the above push down map by
  $$\mathfrak B': \mathfrak B'(\alpha)=\hat{\beta}.$$
 \begin{remark}
  Let us use local coordinates on $U=G_{\mathbb{C}}\times_{H_{\mathbb
  C}}(D_1\times D_2)$ based at a point $x$ with stabilizer group $H$ to see the push down map
  described in the proof of Lemma~\ref{pushdown}. Let $z_0$ be the coordinate along the
 $G_{\mathbb{C}}$-orbit direction, and let $(z_1, z_2)\in D_1\times D_2$ be the coordinate
 in the transversal direction. Then, for instance, we may write a
 $G_{\mathbb{C}}$-invariant $(n, 0)$-form
$\alpha=f(z_0, z_1, z_2)dz_0\wedge dz_1\wedge dz_2$ locally, where
$f$ is a $G_{\mathbb{C}}$-invariant function. Restricting $\alpha$
 to $G_{\mathbb{C}}\cdot (U\cap Z_{(H)})$, we get $\alpha_|=f(z_0, z_1, 0) dz_0\wedge dz_1$.
 The contraction gives $\beta_|= f(z_0, z_1, 0) dz_1$ (up to a sign), and the restriction of $\beta_|$
to $U\cap Z_{(H)}$ gives ${\beta_|}_|= f(x_0, z_1, 0) dz_1$. By the
$G$-invariance of $\alpha$,
 ${\beta_|}_|=\pi^* (\hat{\beta}_|)$, where $\hat{\beta}_|= f([x], z_1, 0) dz_1$ is a local form on $D_1$
 which is a neighborhood of
 $[x]$ in $\mathcal S_{(H)}$. Notice that, conversely, if we have such a $(d_{\mathcal S}, 0)$-form
$\hat{\beta}_|$ on  $\mathcal S_{(H)}$,  we can lift it to a
$(d_{\mathcal S} + d_{G/H}, 0)$-form on
 $G_{\mathbb{C}}\cdot Z_{(H)}$
by using $G_{\mathbb{C}}$-invariance and by ``growing back'' the coordinate $z_0$.
  \end{remark}

  \begin{remark}
  If $G$ acts freely on $\phi^{-1}(0)$, then  $\phi^{-1}(0)=Z_1$, where $1\in G$ is the identity element.
  By the holomorphic slice theorem, a saturated neighborhood of each point
  $x\in \phi^{-1}(0)$ is biholomorphic to $U=G_{\mathbb{C}}\times D$.
  So $U\cap\phi^{-1}(0)=G\times D$, and $D=(U\cap\phi^{-1}(0))/G=U//G_{\mathbb{C}}$ is
  biholomorphic to a neighborhood of $[x]$ in $M_0$. In our point of
  view, we restrict a $G_{\mathbb{C}}$-invariant $(n, 0)$-form $\alpha$ to $U$, then we
  contract the form at the points in $U\cap\phi^{-1}(0)$ with the
  generating vector fields of the
  free $G$-action, then we restrict the resulting form to
  $U\cap\phi^{-1}(0)$ and push it down to $D$ by the quotient map
  $\pi: U\cap\phi^{-1}(0)\rightarrow D$. One may see this in local
  coordinates as we did in the last remark. In the point of view of
  Hall and Kirwin, they contract the form $\alpha$ at the points in
  $U$ with the generating vector fields of the free $G$-action ($G$
  acts freely on $U$), and then push down the resulting form to $D$
  by the quotient map $\pi_{\mathbb C}: U\rightarrow U//G_{\mathbb{C}}=D$. We see that
  the results are the same. Their pulling back of an $(n', 0)$-form
  (let $n'=$dim$(D)$) on $D$ to a  $G_{\mathbb{C}}$-invariant  $(n, 0)$-form to $U$
  is just
  by using the $G_{\mathbb{C}}$-action and by ``adding''  the $G_{\mathbb{C}}$-orbit
  direction.
  \end{remark}

     We define $\hat K$ to be $K//G_{\mathbb{C}}=K|_{\phi^{-1}(0)}/G$, and we define  $\sqrt{\hat
   K}$ to be $\sqrt K//G_{\mathbb{C}}=\sqrt K|_{\phi^{-1}(0)}/G$.
   Sections of $\hat K$ over $M//G$ are ``stratified forms'' $\hat\beta$ whose
   restriction to each stratum $\mathcal S$ of complex dimension
   $d_{\mathcal S}$ is a smooth $(d_{\mathcal S}, 0)$-form.
    If $O$ is a small open set in $M//G$, a section of $\hat K$ over
   $O$ looks like $\hat f_| dw_1\wedge dw_2\wedge...\wedge dw_r$ on
   the  open dense stratum and looks like $\hat f_|
   dw_{i_1}\wedge....\wedge dw_{i_j}$ for some subset $\{i_1, ..., i_j\}$
   of  $\{1, ..., r\}$ on other strata, where $\pi_{\mathbb C}^*(\hat f)$
   is a $G_{\mathbb C}$-invariant function on $\pi_{\mathbb
   C}^*(O)$.\\

   We defined the map
  $\mathfrak B': \Gamma(M, K)^{G_{\mathbb{C}}}\rightarrow   \Gamma(M//G, \hat
   K)$. Using this map, we define a linear map
   $$B': \Gamma(M, \sqrt K)^{G_{\mathbb{C}}}\rightarrow\Gamma(M//G, \sqrt{\hat
     K})$$
 such that
 \begin{equation}\label{B'}
 (B'\mu)^2=\mathfrak B'(\mu^2).
 \end{equation}

  Using the map $A_k'$ and the map $B'$, for each $k$,
 we define a linear map
$$B_k': \Gamma(M, L^{\otimes k}\otimes\sqrt K)^G\rightarrow \Gamma(M//G, (L^{\otimes k})_0\otimes\sqrt{\hat K}),$$
unique up to an overall sign, such that
   $$B_k'(s\otimes\mu)=A_k'(s)\otimes B'(\mu)$$
   where $s\in \Gamma(L^{\otimes k})^G$ and $\mu\in\Gamma(\sqrt K)^G$.\\

   Next, we give an argument of the facts
   $$\mathcal H (M//G, \hat
   K)=\mathfrak B'(\mathcal H (M^{ss}, K)^{G_{\mathbb{C}}}),$$
   and
   $$\mathcal H (M//G, \sqrt{\hat
   K})=B'(\mathcal H (M^{ss}, \sqrt K)^{G_{\mathbb{C}}}).$$

   A holomorphic section
   of $K$ (or of $\sqrt K$) defined over a $G_{\mathbb{C}}$-invariant open set is $G$-invariant if
   and only if it is $G_{\mathbb{C}}$-invariant. Let $\mathcal K$ (or $\sqrt \mathcal K$)
 be the sheaf of holomorphic sections of $K$ (or of $\sqrt K$), and,
   we define a sheaf $\mathcal K'$ (or $\sqrt{\mathcal K'}$)  on $M_0$, by letting
    $\mathcal K'(O)=\mathfrak B'(\mathcal K(\pi_{\mathbb{C}}^{-1}(O))^{G_{\mathbb{C}}})$
  (or by letting
     $\sqrt{\mathcal K'}(O)=B'(\sqrt\mathcal K(\pi_{\mathbb{C}}^{-1}(O))^{G_{\mathbb{C}}})$.) for
    each open set $O$ of $M_0$.
  Using our results above and combining the argument of the proof of Proposition~\ref{propS'}, we have
   the following:

    \begin{proposition}\label{propL}
     The sheaf $\mathcal K'$ (or $\sqrt{\mathcal K'}$) is the sheaf of holomorphic sections of the
     stratified-line bundle $\hat K$ (or $\sqrt{\hat K}$)   over $M_0=M^{ss}//G_{\mathbb{C}}$.
    \end{proposition}

 So we have proved the following
  \begin{theorem}\label{redhfb'}
  Let  $\sqrt{\hat K}=\sqrt K//G_{\mathbb{C}}$.
  There exists a linear map
  $B': \Gamma(M, \sqrt K)^{G_{\mathbb{C}}}\rightarrow \Gamma(M//G, \sqrt{\hat
   K})$, unique up to an overall sign, such that for each
  $\mu\in\Gamma(M, \sqrt K)^{G_{\mathbb{C}}}$, $B'(\mu)$ is a ``stratified form'' on $M//G$
  such that for each stratum $\mathcal S$ of $M//G$ with complex dimension $d_{\mathcal S}$,
  $(B'(\mu))^2|_{\mathcal S}$
  is a $(d_{\mathcal S}, 0)$-
 form on $\mathcal S$, and
 these forms are related by suitable degenerating of directions from higher dimensional
 strata to lower dimensional strata. Moreover, $\mathcal H (M//G, \sqrt{\hat
   K})=B'(\mathcal H (M^{ss}, \sqrt K)^{G_{\mathbb{C}}}).$

  Consequently, for each $k$, there exists a linear map
    $B_k': \Gamma(M, L^{\otimes k}\otimes\sqrt K)^{G_{\mathbb{C}}}\rightarrow \Gamma(M//G, (L^{\otimes k})_0\otimes\sqrt{\hat
   K})$, unique up to an overall sign, such that
   $$B_k'(s\otimes\mu)=A_k'(s)\otimes B'(\mu)$$
   for $s\in \Gamma(M, L^{\otimes k})^{G_{\mathbb{C}}}$ and $\mu\in\Gamma(M, \sqrt K)^{G_{\mathbb{C}}}$,
   and such that $\mathcal H (M//G, (L^{\otimes k})_0\otimes\sqrt{\hat
   K})=B_k'(\mathcal H (M^{ss}, L^{\otimes k}\otimes\sqrt K)^{G_{\mathbb{C}}}).$
  \end{theorem}

  We end this section by giving the definition of a pointwise Hermitian
    structure on $\Gamma(M//G, \sqrt{\hat K})$.
    Let $\mu', \nu'\in\Gamma(M//G, \sqrt{\hat K})$. We define a
    Hermitian
    structure on $\Gamma(M//G, \sqrt{\hat K})$ stratum-wise as
  \begin{equation}\label{norm-red-srb}
    (\mu')^2\wedge(\bar{\nu'})^2|_{\mathcal S}=(\mu',
    \nu')^2|_{\mathcal S}\,\epsilon_{\hat{\omega}_{\mathcal S}},
 \end{equation}
 where $\epsilon_{\hat{\omega}_{\mathcal S}}$ is the volume form on the
    stratum $\mathcal S$ of $M//G$.
  \section{modified linear space isomorphism}

    The following theorem gives the growth of the pointwise norm
    square of a $G$-invariant holomorphic section and a modified
    $G$-invariant holomorphic section along the gradient curves of
    the moment map components. We need this theorem to prove
    Theorem~\ref{modisom}, and we will need this theorem in the
    subsequent sections.

      \begin{theorem}\label{normvary'}
    Let $s\in\mathcal H(M, L^{\otimes k})^G$ and let $r\in\mathcal H(M, L^{\otimes k}\otimes\sqrt K)^G$.
    Let $y_0\in M$, and let $H$ be its stabilizer group. Let $\mathfrak h=Lie(H)$.
    Let $\mathfrak m$ be the orthogonal complement of $\mathfrak h$ in  $\mathfrak g=Lie (G)$
   (assuming we have chosen a $G$-invariant metric on $\mathfrak g$).
    Then for $0\neq\xi\in \mathfrak m$, we have

    (a) $|s|^2(e^{i\xi}\cdot y_0)=|s|^2(y_0)exp\lbrace -\int_0^1
    2k\phi_{\xi}(e^{it\xi}\cdot y_0)dt\rbrace,$

    (b) $|r|^2(e^{i\xi}\cdot y_0)=|r|^2(y_0)exp\lbrace -\int_0^1
    (2k\phi_{\xi}(e^{it\xi}\cdot y_0)+\frac{\mathcal
    L_{JX^{\xi}}\epsilon_{\omega}}{2\epsilon_{\omega}}(e^{it\xi}\cdot
    y_0))dt\rbrace$.

    If we let $f(\xi, y_0):=2\int_0^1\phi_{\xi}(e^{it\xi}\cdot
    y_0)dt$, then as a function of $\xi\in \mathfrak m$, $f(\xi, y_0)$ achieves its
   unique  minimum at $\xi=0$.
   The Hessian of $f(\xi,
     y_0)$ at $\xi=0$ is given by

     $D_{\xi_1}D_{\xi_2}f(\xi, y_0)|_{\xi=0}=2B_{y_0}(JX^{\xi_1},
     JX^{\xi_2})$, $\xi_1, \xi_2\in\mathfrak m$.
    \end{theorem}
    \begin{proof}
    See the proof of Theorem 4.1 in \cite{HK}. Modify the proof by
    noticing the following: for $y_0\in M$, if $H$ is the stabilizer group
    of $y_0$, and if $0\neq\xi \in \mathfrak h=Lie (H)$, then
    $X^{\xi}(y_0)=0$, so  $JX^{\xi}(y_0)=0$ as well, therefore $e^{i\xi}\cdot
    y_0=y_0$.
    \end{proof}

    Using (b) of the above theorem,  we obtain the following modified linear space isomorphism:
    \begin{theorem}\label{modisom}
     For $k$ sufficiently large, the map

     $$B_k': \mathcal H(M, L^{\otimes k}\otimes\sqrt K)^G\rightarrow\mathcal H(M//G,  (L^{\otimes k})_0\otimes\sqrt{\hat
   K})$$

      is bijective.
    \end{theorem}

    \begin{proof}
    We use a similar argument as used by Guillemin and Sternberg in
    \cite{GS}, by Sjamaar in \cite{S}, and by Hall and Kirwin
   in \cite{HK} (the proof of Theorem 3.2 in \cite{HK}).

    By Theorems~\ref{isom} and \ref{redhfb'}, elements in $\mathcal H(M//G,  (L^{\otimes k})_0\otimes\sqrt{\hat
   K})$ lift to elements in $\mathcal H(M^{ss}, L^{\otimes k}\otimes\sqrt
   K)^G$.

     The map is injective because two holomorphic sections which
     agree on the semistable set, which is open and dense in $M$,
     must be equal.

   Let $\hat r\in\mathcal H(M//G,  (L^{\otimes k})_0\otimes\sqrt{\hat
   K})$, and let  $r\in\mathcal H(M^{ss}, L^{\otimes k}\otimes\sqrt
   K)^G$ be its lift. We only need to show that $|r|^2$ remains bounded as we approach
   the unsemistable set (which is of complex codimension at least one), the Riemann Extension Theorem will imply that
   $r$ extends holomorphically to all of $M$.

    By
    Theorem~\ref{normvary'} (b), for  $y_0\in M^{ss}$ with stabilizer group
  $H$, and for $\xi\in\mathfrak m$, we
    have
     $$\frac{d}{dt}|r|^2(e^{it\xi}\cdot y_0)=|r|^2(e^{it\xi}\cdot y_0)
    (-2k\phi_{\xi}(e^{it\xi}\cdot y_0)-\frac{\mathcal L_{JX^{\xi}}\epsilon_{\omega}}{2\epsilon_{\omega}}(e^{it\xi}\cdot y_0)).$$

    Notice that, for $\xi\in\mathfrak h=Lie(H)$, $|r|^2(e^{it\xi}\cdot y_0)=|r|^2(y_0)$, and so $\frac{d}{dt}|r|^2(e^{it\xi}\cdot
    y_0)=0$.

    By compactness of $M$ and by compactness of the set  $\{\xi\in\mathfrak g: |\xi|=1\}$,
    $\frac{L_{JX^{\xi}}\epsilon_{\omega}}{2\epsilon_{\omega}}$ is bounded
    uniformly for all $\xi\in\mathfrak g$ with $|\xi|=1$ and at all points in $M$.

    By the monotonicity of $\phi_{\xi}(e^{it\xi}\cdot y_0)$ in $t$ for $\xi\in\mathfrak
    m$, and by the above fact about $\xi\in\mathfrak h$, we see that
    for all sufficiently
    large $k$,
   $\frac{d}{dt}|r|^2(e^{it\xi}\cdot
    y_0)\leq 0$  for all
    $y_0\in M^{ss}$, all $\xi\in\mathfrak g$ with $|\xi|=1$,
    and all $t\geq 1$.
   It follows that the $r$ obtained extends holomorphically to all of $M$.
    \end{proof}

   \section{the pointwise norms  of the modified sections}

   \begin{theorem}\label{newnorm}
    Suppose $r\in\mathcal H(M, L^{\otimes k}\otimes\sqrt K)^G$.
    Let $x_0\in\phi^{-1}(0)$ be a point with stabilizer group $H$.
    So $[x_0]\in\mathcal S_{(H)}=Z_{(H)}/G$.  Then,

   if $H=G$,  $$\pi^*(|B_k' r|^2([x_0]))=|r|^2(x_0);$$

   otherwise,
    $$\pi^*(|B_k' r|^2([x_0]))=2^{-d_{G/H}/2}vol(G\cdot x_0)|r|^2(x_0),$$
    where $d_{G/H}$ is the dimension of $G/H$.
    \end{theorem}

    By modifying the proofs of Lemma 3.4 and Lemma 3.5 in
    \cite{HK}, we can prove the following two lemmas.
  \begin{lemma}\label{lem1}
   Let $x\in M$ be a point with isotropy group $H$. Let $G\cdot x$
   be the orbit through $x$. Assume that we have chosen a normalized $G$-invariant
  inner product on $\mathfrak g$. Let
  ${\xi_1, ..., \xi_h}$ be an orthonormal basis of $\mathfrak h=Lie(H)$. If $\mathfrak h\neq\mathfrak
  g$, we
  expand ${\xi_1, ..., \xi_h}$ to an orthonormal basis of
 $\mathfrak g$ by joining ${\xi_{h+1}, ..., \xi_d}$. Then, the
 function $\sqrt{det_{j, k=h+1, ..., d}(B(X^{Ad(g)\xi_j}, X^{Ad(g)\xi_k}))_{g\cdot x}}$
 is a constant along the $G$ orbit $G\cdot x$, and
   $$vol(G\cdot x)=\sqrt{det_{j, k=h+1, ..., d}(B(X^{\xi_j}, X^{\xi_k}))_x} .$$
 \end{lemma}
    \begin{lemma}\label{lem2}
    Let $x_0\in\phi^{-1}(0)$. Assume $x_0$ has isotropy group $H\neq G$.
    Choose a basis as in the Lemma above. Let $Z^j=\pi_+X^{\xi_j}=\frac{1}{2}(X^{\xi_j}-iJX^{\xi_j})$, for
    $ j=h+1, ..., d$. Let $\mathcal S_{(H)}=Z_{(H)}/G$, and let
    $d_{\mathcal S_{(H)}}=dim_{\mathbb C} (\mathcal S_{(H)})$. Then $dim_{\mathbb C}
    (G_{\mathbb C}\cdot Z_{(H)})=d_{\mathcal S_{(H)}}+d_{G/H}$.
    Let $\omega_|=\omega|_{G_{\mathbb C}\cdot Z_{(H)}}$, and let
    $(\epsilon_{\omega_|})_|=\frac{\omega_|^{(d_{\mathcal
    S_{(H)}}+d_{G/H})}}{(d_{\mathcal S_{(H)}}+d_{G/H})!}$.
    Then
    $$i(\bigwedge_jZ^j)\circ
    i(\bigwedge_k\bar{Z}^k(\epsilon_{\omega_|})_|(x_0)|_{Z_{(H)}}=2^{-d_{G/H}}vol(G\cdot
    x_0)^2\,\frac{\omega^{d_{\mathcal S_{(H)}}}}{d_{\mathcal
    S_{(H)}}!} (x_0)|_{Z_{(H)}}.$$
    \end{lemma}

    The proof of Lemma~\ref{lem2} uses the result of Lemma~\ref{lem1}.
    Now, we use Lemma~\ref{lem2} to prove Theorem~\ref{newnorm}.

    \begin{proof}
    Near $x_0$, we can write $r=s\mu$, where $s$ is a local
    $G$-invariant holomorphic section of $ L^{\otimes k}$ and $\mu$
    is a local $G$-invariant holomorphic section of $\sqrt K$. Let
    $\alpha=\mu^2$, and let $\alpha_|=\alpha|_{G_{\mathbb C}\cdot Z_{(H)}}$.
   Then, by (\ref{norm-red-srb}) and by (\ref{B'}), we have
   $$(*)\,\,\,\,\,\,\,\,\pi^*((B'\mu, B'\mu)^2\epsilon_{\hat{\omega}_{\mathcal
       S_{(H)}}}([x_0]))=\pi^*(\mathfrak B'(\alpha)([x_0])\wedge\mathfrak
   B'(\bar{\alpha})([x_0])|_{\mathcal S_{(H)}}).$$
   Case 1. Assume $H=G$. Then, by  construction of the map $\mathfrak B'$ (see the proof of Lemma~\ref{pushdown}),
      $$(*)=(\alpha_|\wedge\bar{\alpha_|})(x_0)|_{Z_{G}}\overset{\mbox{by}\,\, (\ref{hermitian-srb})}=(\mu,
      \mu)^2(\epsilon_{\omega_|})_|(x_0)|_{Z_G}=(\mu, \mu)^2\,\frac{\omega^{d_{\mathcal
     S_G}}}{d_{\mathcal S_G}!}(x_0)|_{Z_G}.$$\\
  Case 2. Assume $H\neq G$. Then, by
   construction of the map $\mathfrak B'$,
   $$(*)=i(\bigwedge_{j=h+1, ..., d}Z^j)\alpha_|(x_0)\wedge i(\bigwedge_{k=h+1,
   ..., d}\bar{Z}^k)\bar{\alpha_|}(x_0)|_{Z_{(H)}}.$$
 Since $\alpha$ is holomorphic,
 $i(\pi_+X^{\xi})\bar\alpha=i(\pi_-X^{\xi})\alpha=0$. So the above
$$=(i(\bigwedge_{j=h+1, ..., d}Z^j\wedge\bigwedge_{k=h+1, ...,
  d}\bar{Z}^k)(\alpha_|\wedge\bar{\alpha_|})(x_0))|_{Z_{(H)}}$$
$$\overset{\mbox{by}\,\, (\ref{hermitian-srb})}=i(\bigwedge_{j=h+1, ..., d}Z^j\wedge\bigwedge_{k=h+1, ...,
  d}\bar{Z}^k)((\mu, \mu)^2(\epsilon_{\omega_|})_|)(x_0)|_{Z_{(H)}} $$
 $$=(\mu,\mu)^2(x_0) 2^{-d_{G/H}}vol(G\cdot  x_0)^2\, \frac{\omega^{d_{\mathcal
     S_{(H)}}}}{d_{\mathcal S_{(H)}}!}(x_0)|_{Z_{(H)}}\,\,\mbox{by
     Lemma}~\ref{lem2}.$$
 In the above, we used the fact that, if we do
    $\alpha\wedge\bar{\alpha}=(\mu, \mu)^2\epsilon_{\omega}$ on $M$, we get a
    function $(\mu, \mu)^2$ on $M$, the value $(\mu, \mu)^2(x_0)$ is
    the same as the value  $(\mu,
    \mu)^2_|(x_0)$ of the function $(\mu,
    \mu)^2_|$ obtained by doing   $\alpha_|\wedge\bar{\alpha_|}=(\mu,
    \mu)^2_|(\epsilon_{\omega_|})_|$ on the K\"ahler
    submanifold $G_{\mathbb C}\cdot Z_{(H)}$.\\

  In both cases, dividing by
     $\pi^*\epsilon_{\hat{\omega}_{\mathcal S_{(H)}}}=\omega^{d_{\mathcal S_{(H)}}}/(d_{\mathcal S_{(H)}})!|_{Z_{(H)}}$,
     taking square root and using the fact
     $|A_k's|^2([x_0])=|s|^2(x_0)$, we obtain the result.
     \end{proof}

  \section{Norms of sections in the quantum spaces}
  \subsection{The coarea formula}
   The coarea formula was cited in  \cite{HK}. For convenience, we also include
  this formula here (see \cite{Chav}, pg. 159-160).
  \begin{lemma}
  Let $Q$ and $R$ be smooth Riemannian manifolds with
  dim$(Q)\geq$dim$(R)$, and let $p: Q\rightarrow R$. Then for
  any $f\in L^1(Q)$, one has
  $$\int_Q\mathcal J_pf dvol(Q)=\int_R
  dvol(R)(y)\int_{p^{-1}(y)}(f|_{p^{-1}(y)})dvol(p^{-1}(y))\, ,$$
  where the Jacobian is $\mathcal J_p:=\sqrt{det (p_*\circ
  p_*^{adj})}$.
  \end{lemma}

   For instance, consider $G_{\mathbb C}\cdot Z_{(H)}$ or $G_{\mathbb C}\cdot S_i$
  occuring in
 Lemma~\ref{lem-split1} or in Lemma~\ref{lem-split2} of Section 11.2.
 Let $\mathfrak m$ be the orthogonal complement of $\mathfrak h=Lie(H)$ or of $\mathfrak h'$
in $\mathfrak g$. Denote  $Z_{(H)}$ or  $S_i$ simply as $S$.  Let
$\Lambda: \mathfrak m\times S\rightarrow G_{\mathbb
   C}\cdot S$ be the diffeomorphism $\Lambda(\xi,
   u)=e^{iAd(g)\xi}\cdot u$, if $u$ has infinitesimal isotropy Lie
   algebra  $Ad(g)\mathfrak h$ or
  $Ad(g)\mathfrak h'$.
  The volume element of $G_{\mathbb
   C}\cdot S$ inherited from $M$ at  a point $(\xi,\, g\cdot u)$, where $u$ has isotropy Lie
   algebra $\mathfrak h$ or $\mathfrak h'$, decomposes as
  \begin{equation}\label{eq1}
   {\Lambda}^*(dvol(G_{\mathbb
   C}\cdot S))_{(\xi, \, g\cdot u)}=\tau(\xi, u)dvol(\mathfrak
   m)\wedge dvol(S)_{g\cdot u}
  \end{equation}
   for some $G$-invariant smooth Jacobian function $\tau$, where $dvol(\mathfrak
   m)$ is the Lebesgue measure on $\mathfrak m$.

   \subsection{Norms of sections in the quantum spaces}

  In this section, we compute the norms of the sections in the quantum spaces.
  The main result of this section is Theorem~\ref{qnsr}.\\

   Let $Z_{(H)}$ be as in (\ref{Z_H}).
 Then
   $$M^{ss}=\bigcup_{(H)} F_{\infty}^{-1} (Z_{(H)}),$$
where $F_{\infty}$ is the limit map of
 the flow $F_t$
 of the gradient of $-\|\phi\|^2$. By dividing further into connected components for each
   $Z_{(H)}$, we assume
   that each $Z_{(H)}$ is connected.     Since there is an open dense connected
   stratum  $Z_{(H)}$ (for some $H$) in $\phi^{-1}(0)$ (\cite{SL}), there is an open dense
  connected set
  $F_{\infty}^{-1} (Z_{(H)})$ in  $M^{ss}$.
 We will compute the
  integral of the pointwise norm square of the
sections over each $F_{\infty}^{-1}
  (Z_{(H)})$ with $H$ varing. This integral will relate to the integral over
   $\mathcal S_{(H)}$ (see (\ref{S_H}) for the definition of $\mathcal S_{(H)}$)      of the
  pointwise norm square of the descended section. In particular, if $Z_{(H)}$ or
 $\mathcal S_{(H)}$ is a single point, then the integral on it is regarded as
 the pointwise norm square of the sections over this point.\\

   Now, let us take a look at  $F_{\infty}^{-1}(Z_{(H)})$. By 2. and 3. of Proposition~\ref{propS},
   $G_{\mathbb{C}}\cdot Z_{(H)}\subseteqq F_{\infty}^{-1}(Z_{(H)})$.
  By the holomorphic slice theorem or by
  Theorem~\ref{theorem-stable}, we have
\begin{lemma}
   If $d\phi_x$ is
  surjective for all $x\in Z_{(H)}$ ($H$ is necessarily finite), then
$F_{\infty}^{-1}(Z_{(H)})=G_{\mathbb{C}}\cdot Z_{(H)}$.
\end{lemma}
Now, we assume that
 $F_{\infty}^{-1}(Z_{(H)})$ is strictly
   larger than $G_{\mathbb{C}}\cdot Z_{(H)}$.
 Then $F_{\infty}^{-1}(Z_{(H)})$ contains
  complex  orbits whose closures contain the complex orbits in $G_{\mathbb{C}}\cdot Z_{(H)}$.
We decompose the set $F_{\infty}^{-1}(Z_{(H)})-G_{\mathbb{C}}\cdot
Z_{(H)}$ into a
 disjoint union of
 connected $G_{\mathbb{C}}$-invariant complex submanifolds each of which has
 the same infinitesimal compact orbit type, say $\mathfrak{(h')}$. Let
   $M^{(H)}_{\mathfrak{(h')}}$ denote one of these invariant complex submanifold.

\begin{lemma}\label{lem-dimension}
Assume that  $F_{\infty}^{-1}(Z_{(H)})$ is strictly
   larger than $G_{\mathbb{C}}\cdot Z_{(H)}$. We decompose
$F_{\infty}^{-1}(Z_{(H)})-G_{\mathbb{C}}\cdot
Z_{(H)}=\bigcup_{\mathfrak  h'}M^{(H)}_{\mathfrak{(h')}}$, where
$\bigcup_{\mathfrak  h'}M^{(H)}_{\mathfrak{(h')}}$ is a disjoint
union with each  $M^{(H)}_{\mathfrak{(h')}}$ being a connected
$G_{\mathbb{C}}$-invariant complex submanifold of a certain
infinitesimal compact orbit type $\mathfrak{(h')}$.
 Then $0\notin\phi(M^{(H)}_{\mathfrak{(h')}})$ and  dim $(\mathfrak h')< dim (H)$. Therefore, for any
possible $H$, if $\mathfrak m$ is the orthogonal complement of
$\mathfrak h'$ in $\mathfrak g$, then dim $(\mathfrak m)>0$.
\end{lemma}

\begin{proof}
   By 2. of Proposition~\ref{propS}, the complex orbits in  $M^{(H)}_{\mathfrak{(h')}}$
    do not intersect $\phi^{-1}(0)$.
 So $0\notin\phi(M^{(H)}_{\mathfrak{(h')}})$.

 By Theorem~\ref{slice}, a neighborhood $U$ of $x\in Z_{(H)}$ in
  $M$ is $G$-equivariantly biholomorphic to $G_{\mathbb{C}}\times_{H_{\mathbb{C}}} D$.
  Split $D=D_1\times D_2$,
   where $D_1$ is the complex subspace fixed by $H$ and $H_{\mathbb{C}}$.
 Let   $\phi_|$ be the
    moment map of the $H$ action on $D_2$ with respect to the restricted K\"ahler form.
    Then $U\cap\phi^{-1}(0)=G\times_H (D_1\times \phi_|^{-1}(0))$,
  $U\cap Z_{(H)}=G\times_H (D_1\times
 0)$, and $U\cap G_{\mathbb{C}}\cdot Z_{(H)}=G_{\mathbb{C}}\times_{H_{\mathbb{C}}}
 (D_1\times 0)$. By the assumption, $D_2\neq \emptyset$ and $M^{(H)}_{\mathfrak{(h')}}\cap
 D_2\neq\emptyset$. The set $M^{(H)}_{\mathfrak{(h')}}\cap
 D_2$ is an $H$-invariant subset of $D_2$ consisting of points with isotropy Lie algebra
 $\mathfrak h'\subset \mathfrak h=Lie(H)$ (A group $H'$ such that $Lie(H')=\mathfrak h'$ is
 a subgroup of $H$ since any point in $D_2$ has isotropy group a subgroup of $H$).

 Since $\phi(M^{(H)}_{\mathfrak{(h')}})$ does not intersect $0$,
 $\phi_|(M^{(H)}_{\mathfrak{(h')}}\cap D_2)$ does not intersect $0$.
 One only needs to argue when $H$ is not connected and when dim $(H')$=dim $(H)$ and
 exclude this possibility by using the fact that
 a finite group action does not
  contribute to the moment map $\phi_|$.

\end{proof}

      By definition of the moment map,
 for $x\in M$ with isotropy Lie algebra $\mathfrak h'$, the image of $d\phi_x: T_xM \rightarrow\mathfrak{g}^*$
       is the annihilator
       in $\mathfrak{g}^*$ of $\mathfrak h'$. So the image $\phi(M^{(H)}_{\mathfrak{(h')}})$
       intersects with a closed positive Weyl chamber at a certain dimension. This image may lie on one or more than one open faces of the moment polytope
 $\triangle$ of $\phi$. These faces form a connected set since we took $M^{(H)}_{\mathfrak{(h')}}$
 to be connected. For a non-abelian Lie group action, the moment
 polytope is defined to be the intersection of the image of the
 moment map with a fixed closed positive Weyl chamber. The faces of
 the moment polytope are caused by symplectic submanifolds with
 different isotropy groups. One should distinguish the faces of the
 moment polytope with the faces of the Weyl chamber.

\begin{lemma}\label{lem-split1}
  Assume $\phi(M^{(H)}_{\mathfrak{(h')}})$ only lies on one open face  $\mathcal
  F_0$ of
  $\triangle$. Then dim $(\mathcal F_0)>0$. Let $0\neq a_0\in\phi(M^{(H)}_{\mathfrak{(h')}})\subset \mathcal F_0$ be
  a value.
 Then, we can write $M^{(H)}_{\mathfrak{(h')}}=G_{\mathbb C}\cdot S_0$, where
 $S_0\subseteqq S_{\mathfrak{(h')}}=\{x\in\phi^{-1}(G\cdot a_0): x$ has isotropy Lie algebra type $(\mathfrak
 h')\}$ and $S_0$ is $G$-invariant.
\end{lemma}

\begin{proof}
 Since  $F_{\infty}(M^{(H)}_{\mathfrak{(h')}})\subset\phi^{-1}(0)$,
there are points
 in $\phi(M^{(H)}_{\mathfrak{(h')}})$ arbitrarily near $0$. Since
 $\phi(M^{(H)}_{\mathfrak{(h')}})$ is connected, dim $(\mathcal F_0)>0$.

 Since $S_0\subset M^{(H)}_{\mathfrak{(h')}}$, and since $M^{(H)}_{\mathfrak{(h')}}$ is
 $G_{\mathbb C}$-invariant, we have $G_{\mathbb C}\cdot S_0\subset M^{(H)}_{\mathfrak{(h')}}$.
Conversely, if $x\in M^{(H)}_{\mathfrak{(h')}}$, then $\phi(x)\in
G\cdot\mathcal F_0$. Without loss of generality, we assume the
isotropy Lie algebra of $x$ is $\mathfrak h'$ and
$\phi(x)=b\in\mathcal F_0$.  If $b=a_0$, then $x\in
S_{\mathfrak{(h')}}$.
 If $b\neq a_0$, then $x$ can be reached by the flow line of $JX^{\xi}$ from a point in $\phi^{-1}(a_0)$, where
 $\xi$ is a vector in $\mathcal F_0$ pointing from $a_0$ to $b$. So $x\in G_{\mathbb C}\cdot S_{\mathfrak{(h')}}$.
\end{proof}

\begin{corollary}
 Let $a$ be any point on the face $\mathcal F_0$, and let
 $S\subseteqq \{x\in\phi^{-1}(G\cdot a): x$ has isotropy Lie algebra type $(\mathfrak
 h')\}$. Then 
 $G_{\mathbb C}\cdot S_0=G_{\mathbb C}\cdot S=M^{(H)}_{\mathfrak{(h')}}$.
\end{corollary}
\begin{proof}
  We have $\mathcal F_0\subset\mathfrak m$,
 where $\mathfrak m$ is the orthogonal complement of
$\mathfrak h'$ in $\mathfrak g$ which is identified with the
annihilator of $\mathfrak h'$  in $\mathfrak{g}^*$.
 The image $\phi(G_{\mathbb C}\cdot S_0)$ must cover the face $\mathcal F_0$.
 So $a\in\phi(M^{(H)}_{\mathfrak{(h')}})$. 
 \end{proof}

 \begin{lemma}\label{lem-split2}
   Assume $\phi(M^{(H)}_{\mathfrak{(h')}})$ lies on
   more than one faces  of
   $\triangle$. Let $\mathcal F_1, ..., \mathcal F_p$ be the ones
   whose closures contain $0$.
  Let $0\neq a_i\in\mathcal F_i$, $i=1, ..., p$, and let
  $$S_i\subseteqq\{x\in\phi^{-1}(G\cdot a_i): x \,\,\mbox{has isotropy Lie algebra type} (\mathfrak
  h')\}.$$
  If $\mathcal F_k$ is in the closure of $\mathcal F_i$, then
$G_{\mathbb C}\cdot S_k\subset G_{\mathbb C}\cdot S_i$. Moreover, we
can write
  $M^{(H)}_{\mathfrak{(h')}}=\cup_{i\in I} (G_{\mathbb C}\cdot
  S_i)$, where $I$ is the subset of $\{1,...p\}$ such that  $\mathcal F_{i\in I}$
  are the top dimensional faces among the $\mathcal F_i$'$s$.
\end{lemma}

 \begin{proof}
 We have $\mathcal F_k\subset \bar{\mathcal F_i}\subset\mathfrak m$,
 where $\mathfrak m$ is the orthogonal complement of
$\mathfrak h'$ in $\mathfrak g$ which is identified with the
annihilator of $\mathfrak h'$  in $\mathfrak{g}^*$.
 Since the points in $S_k$ have isotropy Lie algebra
 $\mathfrak{(h')}$, the $G_{\mathbb C}$ action (or the
 $i\mathfrak{(m)}$ action) will take the points in $S_k$ out
 and emerge them into $G_{\mathbb C}\cdot S_i$. Or, equivalently,
 the moment map value increases  along the flow lines of $JX^{\xi}$,
 where $\xi\in\mathfrak m$ is orthogonal to $\mathcal F_k$. This
 proves $G_{\mathbb C}\cdot S_k\subset G_{\mathbb C}\cdot S_i$.

 So, using Lemma~\ref{lem-split1}, $\cup_{i\in I} (G_{\mathbb C}\cdot
  S_i)=\cup_{i=1}^p (G_{\mathbb C}\cdot
  S_i)\subset M^{(H)}_{\mathfrak{(h')}} $. If $\phi(M^{(H)}_{\mathfrak{(h')}})$
  lies on another face $\mathcal F_{p+1}$ whose closure does not contain $0$, then
   $G_{\mathbb C}\cdot S_{p+1}$ (where $S_{p+1}$ is taken similarly as the $S_i$'$s$)
   should emerge into   $\cup_{i=1}^p G_{\mathbb C}\cdot
   S_i$  to converge to $\phi^{-1}(0)$. This proves
   $M^{(H)}_{\mathfrak{(h')}}=\cup_{i\in I} (G_{\mathbb C}\cdot
   S_i)$.

 \end{proof}

By this lemma, if two faces $\mathcal F_i$ and $\mathcal F_j$ where
$i, j\in I$ contain a one dimensional less face $\mathcal F_k$ in
their common
 closure, then $G_{\mathbb C}\cdot S_i\cap G_{\mathbb C}\cdot S_j=G_{\mathbb C}\cdot S_k$.

 \begin{remark}
 In the above lemma, generally we cannot get all $G_{\mathbb C}\cdot S_i$ from
  $G_{\mathbb C}\cdot S_k$ by the flow lines of $JX^{\xi}$,
 where $\xi\in\mathfrak (m)$. Some orbits in $G_{\mathbb C}\cdot
 S_i$ may converge to more singular orbits in
 $\phi^{-1}(G\cdot\mathcal F_k)$.
 \end{remark}

 Since $0$ is in the closure of each $\mathcal F_i$, dim$(\mathcal F_i)>0$ for each $i=1, ...,
 p$.\\

 So we have proved
 \begin{lemma}\label{D_2}
 We can decompose $F_{\infty}^{-1}(Z_{(H)})$  into a (finite) disjoint union
 $F_{\infty}^{-1}(Z_{(H)})=G_{\mathbb{C}}\cdot Z_{(H)}\bigcup_{\mathfrak  h'}M^{(H)}_{\mathfrak{(h')}}$,
  where $\bigcup_{\mathfrak  h'}M^{(H)}_{\mathfrak{(h')}}=\emptyset$, or, each $M^{(H)}_{\mathfrak{(h')}}$ can be written as in Lemma~\ref{lem-split1} or in Lemma~\ref{lem-split2}.
 In the second case, for any $i=0, 1...,p$, we may choose  $a'_i\neq 0$ on $\mathcal F_i$ different from $a_i$ and
  choose  $S'_i\subseteqq\{x\in\phi^{-1}(G\cdot a'_i): x$ has isotropy Lie algebra type $(\mathfrak h')\}$
 and we have $G_{\mathbb C}\cdot S_i=G_{\mathbb C}\cdot S'_i$.
 \end{lemma}

 \begin{definition}\label{defI-II}
Let $n_{(H)}$ be the complex dimension of
 $G_{\mathbb{C}}\cdot  Z_{(H)}$, and let $n^{(H)}_{\mathfrak{(h')}}$ be the complex dimension of
 $M^{(H)}_{\mathfrak{(h')}}$.
 Take $s\in\mathcal H(M, L^{\otimes k})^G$. Define
   $$I_k^{Z_{(H)}}=(k/2\pi)^{n_{(H)}/2}\int_{G_{\mathbb{C}}\cdot  Z_{(H)}} |s|^2  dvol(G_{\mathbb{C}}\cdot  Z_{(H)}),$$
   and,
$$II_k^{Z_{(H)}}=\sum_{\mathfrak  h'}(k/2\pi)^{n^{(H)}_{\mathfrak{(h')}}/2}\int_{M^{(H)}_{\mathfrak{(h')}}}|s|^2 dvol(M^{(H)}_{\mathfrak{(h')}}).$$
Define
 $$\int_{F_{\infty}^{-1}(Z_{(H)})}|s|^2 dvol(F_{\infty}^{-1}(Z_{(H)}))$$
$$=I_k^{Z_{(H)}}+II_k^{Z_{(H)}}.$$
 \end{definition}

   \begin{lemma}\label{intstr}
 Let $s\in\mathcal H(M, L^{\otimes k})^G$.
 Then\\

(a). $I_k^{Z_{(H)}}=(k/2\pi)^{d_{\mathcal S_{(H)}}/2}\int_{\mathcal
S_{(H)}}|A_k's|^2([x])I_k^{\mathcal
S_{(H)}}([x])\epsilon_{\hat{\omega}_{\mathcal S_{(H)}}},$ where
 $$I_k^{\mathcal
  S_{(H)}}([x])=1,\,\, \mbox{if}\,\, H=G;\,\, \mbox{and}$$
    $$I_k^{\mathcal S_{(H)}}([x])=vol(G\cdot x)(k/2\pi)^{d_{G/H}/2}\int_{\mathfrak m}
   \tau(\xi, x)exp\lbrace -\int_0^1
    2k\phi_{\xi}(e^{it\xi}\cdot x)dt\rbrace
   dvol(\mathfrak m),$$ if $H\neq G$. Here, $\mathfrak m$
   denotes the orthogonal complement of $\mathfrak h=Lie(H)$ in $\mathfrak
   g$, and $x$ is a point with stabilizer group $H$.

 (b). $$II_k^{Z_{(H)}}=0,$$ or
 $$II_k^{Z_{(H)}}=\sum_{\mathfrak  h'}(k/2\pi)^{n^{(H)}_{\mathfrak{(h')}}/2}(\sum_i \pm \int_{S_i}|s|^2 (g\cdot u) dvol(S_i)$$
     $$\int_{\mathfrak m'}\tau(\zeta, u)exp\lbrace -2k\int_0^1\phi_{\zeta}(e^{it\zeta}\cdot u)\rbrace dvol(\mathfrak
   m')),$$
 where the second sum is over some subset of indices of $i$ occuring in Lemma~\ref{lem-split1} or
 in Lemma~\ref{lem-split2}, $\mathfrak m'$
   is the orthogonal complement of $\mathfrak h'$ in $\mathfrak
   g$, and the points $u\in S_i$ are of isotropy Lie algebra $\mathfrak h'$.
\end{lemma}

   \begin{proof}

 We will drop the subscripts and superscripts in $I_k^{Z_{(H)}}$ and $II_k^{Z_{(H)}}$
and simply write $I$ and $II$.\\

(a).    If $H=G$, then  $G_{\mathbb{C}}\cdot Z_{G}=Z_{G}=\mathcal
S_G$. So
   $$I=(k/2\pi)^{n_{(H)}/2}\int_{G_{\mathbb{C}}\cdot
      Z_{G}}|s|^2 dvol(G_{\mathbb{C}}\cdot
      Z_{G})=(k/2\pi)^{n_{(H)}/2}\int_{Z_{G}}|s|^2(x)dvol({Z_{G}})$$
 $$=(k/2\pi)^{d_{\mathcal S_G}/2}\int_{\mathcal S_G}|A_k's|^2([x])\epsilon_{\hat{\omega}_{\mathcal S_G}}.$$

If $H\neq G$, by the coarea formula, the formula (\ref{eq1}), and
Theorem~\ref{normvary'} (a),  we have

 $$I=(k/2\pi)^{n_{(H)}/2}\int_{Z_{(H)}}|s|^2(x')dvol(Z_{(H)})$$
 $$\int_{\mathfrak m}\tau(\xi, g^{-1}x')exp\lbrace
   -2k\int_0^1\phi_{Ad(g)\xi}(e^{itAd(g)\xi}\cdot x')\rbrace
   dvol(\mathfrak m),$$
   where $x'$ is any point in $Z_{(H)}$ with stabilizer group
   $gHg^{-1}$ (for some $g$).

  By $G$-invariance of the function $\tau$, and by $G$-equivariance
  of the moment map $\phi$, we have
  $$\int_{\mathfrak m}\tau(\xi, g^{-1}x')exp\lbrace -2k\int_0^1\phi_{Ad(g)\xi}(e^{itAd(g)\xi}\cdot x')\rbrace
   dvol(\mathfrak m)$$
   $$=\int_{\mathfrak m}\tau(\xi, x)exp\lbrace
   -2k\int_0^1\phi_{\xi}(e^{it\xi}\cdot x)\rbrace
   dvol(\mathfrak m),$$
    where $x=g^{-1}x'$ has stabilizer $H$.

 Using the fact that $dvol(Z_{(H)})=dvol(G\cdot x)\wedge \pi^*dvol(\mathcal S_{(H)})$,
 the integral
$$I=(k/2\pi)^{(n_{(H)}-d_{G/H})/2}\int_{\mathcal S_{(H)}}|A_k's|^2([x])dvol(\mathcal S_{(H)})vol(G\cdot
x)(k/2\pi)^{d_{G/H}/2}$$
$$\int_{\mathfrak m} \tau(\xi, x)exp\lbrace  -2k\int_{\gamma_{\xi}}\phi_{\xi}\rbrace dvol(\mathfrak m)$$
$$=(k/2\pi)^{d_{\mathcal S_{(H)}}/2}\int_{\mathcal S_{(H)}}|A_k's|^2([x])I_k^{\mathcal S_{(H)}}([x])\epsilon_{\hat{\omega}_{\mathcal S_{(H)}}},$$

where    $$I_k^{\mathcal S_{(H)}}([x])=vol(G\cdot
x)(k/2\pi)^{d_{G/H}/2}\int_{\mathfrak m}
   \tau(\xi, x)exp\lbrace
   -2k\int_0^1\phi_{\xi}(e^{it\xi}\cdot x)\rbrace
   dvol(\mathfrak m),$$ with $x$ being taken as a point (on the orbit $G\cdot x$) with stabilizer group
   exactly $H$.\\

(b). If  $\bigcup_{\mathfrak
h'}M^{(H)}_{\mathfrak{(h')}}=\emptyset$, then $II=0$.

Otherwise, let
us only consider one summand for the first summation in $II$. The
others follow similarly.
   So, we assume
$$II=(k/2\pi)^{n^{(H)}_{\mathfrak{(h')}}/2}\int_{M^{(H)}_{\mathfrak{(h')}}}|s|^2
dvol(M^{(H)}_{\mathfrak{(h')}}).$$
 By Lemma~\ref{lem-split1} or Lemma~\ref{lem-split2}, we can compute
 this integral over one set $G_{\mathbb C}\cdot S_0$, or we can
 compute it over a finite union $G_{\mathbb C}\cdot S_{i\in I}$ and possibly subtract
 some
 integrals over some mutual intersections which have  similar forms
 (if a mutual intersection has less dimension, then we do not subtract). So
 we only need to write one such integral in the stated form.

 Using the coarea
formula, the formula (\ref{eq1}) and Theorem~\ref{normvary'} (a) on
the space $G_{\mathbb C}\cdot S_i$, we have
$$\int_{G_{\mathbb C}\cdot S_i}|s|^2dvol(G_{\mathbb C}\cdot S_i)=\int_{S_i}|s|^2 (u') dvol(S_i)$$
$$\int_{\mathfrak m'}
 \tau(\zeta, g^{-1}u')exp\lbrace
   -2k\int_0^1\phi_{Ad(g)\zeta}(e^{itAd(g)\zeta}\cdot u')\rbrace dvol(\mathfrak m'),$$
   where $u'\in S_i$ is any point with isotropy Lie algebra
   $Ad(g)\mathfrak h'$ (for some $g$).

  For the same reason as in (a), we have
 $$\int_{\mathfrak m'}
 \tau(\zeta, g^{-1}u')exp\lbrace
   -2k\int_0^1\phi_{Ad(g)\zeta}(e^{itAd(g)\zeta}\cdot u')\rbrace dvol(\mathfrak m')$$
 $$=\int_{\mathfrak m'}
 \tau(\zeta, u)exp\lbrace
   -2k\int_0^1\phi_{\zeta}(e^{it\zeta}\cdot u)\rbrace dvol(\mathfrak m'),$$
   where $u=g^{-1}u'$ has isotropy Lie algebra $\mathfrak h'$.
 \end{proof}

 \begin{definition}\label{defI'-II'}
 We use the same notations as those in Definition~\ref{defI-II}.
 Take $r\in\mathcal H(M, L^{\otimes k}\otimes\sqrt K)^G$. Define
   $$\tilde I_k^{Z_{(H)}}=(k/2\pi)^{n_{(H)}/2}\int_{G_{\mathbb{C}}\cdot  Z_{(H)}} |r|^2  dvol(G_{\mathbb{C}}\cdot  Z_{(H)}),$$
   and,
 $$\widetilde{II}_k^{Z_{(H)}}=\sum_{\mathfrak  h'}(k/2\pi)^{n^{(H)}_{\mathfrak{(h')}}/2}\int_{M^{(H)}_{\mathfrak{(h')}}}|r|^2 dvol(M^{(H)}_{\mathfrak{(h')}}).$$
 Define
 $$\int_{F_{\infty}^{-1}(Z_{(H)})}|r|^2 dvol(F_{\infty}^{-1}(Z_{(H)}))$$
$$=\tilde I_k^{Z_{(H)}}+\widetilde{II}_k^{Z_{(H)}}.$$
 \end{definition}

 \begin{lemma}\label{intstr'}
Let $r\in\mathcal H(M, L^{\otimes k}\otimes\sqrt K)^G$.  Then\\

(a). $\tilde I_k^{Z_{(H)}}=(k/2\pi)^{d_{\mathcal
S_{(H)}}/2}\int_{\mathcal S_{(H)}}|B_k'r|^2([x])J_k^{\mathcal
S_{(H)}}([x])\epsilon_{\hat{\omega}_{\mathcal S_{(H)}}},$
 where
  $$J_k^{\mathcal
  S_{(H)}}([x])=1,\,\, \mbox{if}\,\, H=G;\,\, \mbox{and}$$
$$J_k^{\mathcal S_{(H)}}([x])=(k/2\pi)^{d_{G/H}/2}2^{d_{G/H}/2}\cdot$$
$$\int_{\mathfrak m}
   \tau(\xi, x)exp\lbrace
   -\int_0^1(2k\phi_{\xi}(e^{it\xi}\cdot x)+\frac{\mathcal
    L_{JX^{\xi}}\epsilon_{\omega}}{2\epsilon_{\omega}}(e^{it\xi}\cdot
    x))\rbrace dvol(\mathfrak m),$$ if $H\neq G$. Here, $\mathfrak m$
   denotes the orthogonal complement of $\mathfrak h=Lie(H)$ in $\mathfrak
   g$, and $x$ is a point with stabilizer group $H$.

 (b). $$\widetilde{II}_k^{Z_{(H)}}=0,$$ or
 $$\widetilde{II}_k^{Z_{(H)}}=\sum_{\mathfrak  h'}(k/2\pi)^{n^{(H)}_{\mathfrak{(h')}}/2}(\sum_i \pm \int_{S_i}|r|^2 (g\cdot u) dvol(S_i)$$
     $$\int_{\mathfrak m'}\tau(\zeta, u)exp\lbrace -\int_0^1(2k\phi_{\zeta}(e^{it\zeta}\cdot u)+\frac{\mathcal L_{JX^{\zeta}}\epsilon_{\omega}}{2\epsilon_{\omega}}(e^{it\zeta}\cdot u))\rbrace
dvol(\mathfrak m')),$$
 where the second sum is over some subset of indices of $i$ occuring in Lemma~\ref{lem-split1}
 or Lemma~\ref{lem-split2}, $\mathfrak m'$
   is the orthogonal complement of $\mathfrak h'$ in $\mathfrak
   g$, and the points $u\in S_i$ are of isotropy Lie algebra $\mathfrak h'$.

\end{lemma}

  \begin{proof}
 (a). The proof is  similar to the proof of (a) of Lemma~\ref{intstr}, but we will use
 Theorem~\ref{newnorm}. We will drop the subscript and superscript in $\tilde I_k^{Z_{(H)}}$  and simply write
 $\tilde I$.\\

    If $H=G$, then  $G_{\mathbb{C}}\cdot
      Z_{G}=Z_{G}=\mathcal S_{G}$. Then

          $$\tilde I=(k/2\pi)^{n_{G}/2}\int_{G_{\mathbb{C}}\cdot Z_{G}}|r|^2
          dvol(G_{\mathbb{C}}\cdot Z_{G})$$
$$=(k/2\pi)^{n_{G}/2}\int_{Z_{G}}|r|^2(x)dvol({Z_{G}})
=(k/2\pi)^{d_{\mathcal S_{G}}/2}\int_{\mathcal
S_{G}}|B_k'r|^2([x])\epsilon_{\hat{\omega}_{\mathcal S_{G}}}$$ by
Theorem~\ref{newnorm}.\\

 If $H\neq G$,  by the coarea
formula, the formula (\ref{eq1}), Theorem~\ref{normvary'} (b), and
by a $G$-invariance argument as in the proof of Lemma~\ref{intstr},
we have

 $$\tilde I=(k/2\pi)^{n_{(H)}/2}\int_{Z_{(H)}}|r|^2(g\cdot x)dvol(Z_{(H)})$$
 $$\int_{\mathfrak m}\tau(\xi, x)exp\lbrace
   -\int_{\gamma_{\xi}}(2k\phi_{\xi}+\frac{\mathcal
    L_{JX^{\xi}}\epsilon_{\omega}}{2\epsilon_{\omega}})\rbrace
   dvol(\mathfrak m).$$
   By  Theorem~\ref{newnorm},
$$\tilde I=(k/2\pi)^{d_{\mathcal S_{(H)}}/2}\int_{\mathcal S_{(H)}}|B_k'r|^2([x])dvol(\mathcal S_{(H)})2^{d_{G/H}/2}(k/2\pi)^{d_{G/H}/2}$$
$$\int_{\mathfrak m}
   \tau(\xi, x)exp\lbrace
   -\int_{\gamma_{\xi}}(2k\phi_{\xi}+\frac{\mathcal
    L_{JX^{\xi}}\epsilon_{\omega}}{2\epsilon_{\omega}})\rbrace
   dvol(\mathfrak m)$$
$$=(k/2\pi)^{d_{\mathcal S_{(H)}}/2}\int_{\mathcal S_{(H)}}|B_k'r|^2([x])J_k^{\mathcal S_{(H)}}([x])\epsilon_{\hat{\omega}_{\mathcal S_{(H)}}},$$

where $$J_k^{\mathcal S_{(H)}}([x])=
(k/2\pi)^{d_{G/H}/2}2^{d_{G/H}/2}\int_{\mathfrak m}
   \tau(\xi, x)exp\lbrace
   -\int_{\gamma_{\xi}}(2k\phi_{\xi}+\frac{\mathcal
    L_{JX^{\xi}}\epsilon_{\omega}}{2\epsilon_{\omega}})\rbrace dvol(\mathfrak m).$$\\

(b). Similar to the proof of (b) of  Lemma~\ref{intstr}. We omit it.
\end{proof}

  Now, we come to our main result of this section:

   \begin{theorem}\label{qnsr}
(a). Let $s\in\mathcal H(M, L^{\otimes k})^G$. Then
 $$\sum_{Z_{(H)}}\int_{F_{\infty}^{-1} (Z_{(H)})}|s|^2 dvol(F_{\infty}^{-1} (Z_{(H)}))$$
 $$=\sum_{\mathcal S_{(H)}}
     (k/2\pi)^{d_{\mathcal S_{(H)}}/2}\int_{\mathcal
       S_{(H)}}|A_k's|^2([x])I_k^{\mathcal
       S_{(H)}}([x])\epsilon_{\hat{\omega}_{\mathcal S_{(H)}}}+\sum_{Z_{(H)}} II_k^{Z_{(H)}},$$
   where   $I_k^{\mathcal S_{(H)}}([x])$ is as in Lemma~\ref{intstr} (a), and each
   $II_k^{Z_{(H)}}$ is as in Lemma~\ref{intstr} (b).\\

 In particular, the above
 is true for each individual summand with respect to $(H)$.\\

(b). Let $r\in\mathcal H(M, L^{\otimes k}\otimes\sqrt K)^G$.  Then
  $$\sum_{Z_{(H)}}\int_{F_{\infty}^{-1} (Z_{(H)})}|r|^2 dvol(F_{\infty}^{-1}(Z_{(H)}))$$
  $$=\sum_{\mathcal S_{(H)}}(k/2\pi)^{d_{\mathcal S_{(H)}}/2}\int_{\mathcal
    S_{(H)}}|B_k'r|^2([x])J_k^{\mathcal
    S_{(H)}}([x])\epsilon_{\hat{\omega}_{\mathcal S_{(H)}}}+\sum_{Z_{(H)}} \widetilde{II}_k^{Z_{(H)}},$$
 where   $J_k^{\mathcal S_{(H)}}([x])$  is as in Lemma~\ref{intstr'} (a), and each
   $\widetilde{II}_k^{Z_{(H)}}$ is as in Lemma~\ref{intstr'} (b).\\

  In particular, the above
 is true for each individual summand with respect to $(H)$.
 \end{theorem}

 \begin{proof}
  Lemmas~\ref{intstr} and \ref{intstr'} proved the statements for the individual
  summands.
 The statements for the sums follow from these lemmas by taking the sum of the individual terms.
 \end{proof}

 The asymptotic properties of $I_k^{\mathcal
  S_{(H)}}([x])$, of $J_k^{\mathcal S_{(H)}}([x])$, of $II_k^{Z_{(H)}}$ and of      $\widetilde{II}_k^{Z_{(H)}}$ will be studied in the
next section (see  Theorem~\ref{asym}).

\section{Asymptotics}

   Our main result of this section is
   \begin{theorem}\label{asym}
  (a).   The densities $I_k^{\mathcal S_{(H)}}$ and $J_k^{\mathcal S_{(H)}}$ for $H\neq G$ satisfy

   $$\mbox{lim}_{k\rightarrow\infty}I_k^{\mathcal S_{(H)}}([x])=2^{-d_{G/H}/2}vol(G\cdot
    x),$$  and
   $$\mbox{lim}_{k\rightarrow\infty}J_k^{\mathcal S_{(H)}}([x])=1.$$
    The limits are uniform for $[x]\in Z_{(H)}/G$.\\

  (b). If $II_k^{Z_{(H)}}\neq 0$ and $\widetilde{II}_k^{Z_{(H)}}\neq 0$, then they
       satisfy
    $$\mbox{lim}_{k\rightarrow\infty}\,II_k^{Z_{(H)}}=0,$$ and
    $$\mbox{lim}_{k\rightarrow\infty}\,\widetilde{II}_k^{Z_{(H)}}=0.$$

   \end{theorem}

   The proof of this theorem will be in Section 12.3.

  \subsection{Growth estimates}
 In Lemmas~\ref{intstr} and \ref{intstr'},
 in the expressions of $I_k^{\mathcal S_{(H)}}$, or of $J_k^{\mathcal
 S_{(H)}}$ (for $H\neq G$), or in the summands of $II_k^{Z_{(H)}}$ or $\widetilde{II}_k^{Z_{(H)}}$, we
 had the following types of integrals

  $$\int_{\mathfrak m}\tau(\xi, x)exp\lbrace
   -2k\int_{\gamma_{\xi}}\phi_{\xi}\rbrace dvol(\mathfrak m)$$
 and
  $$\int_{\mathfrak m}
   \tau(\xi, x)exp\lbrace
   -\int_{\gamma_{\xi}}(2k\phi_{\xi}+\frac{\mathcal
    L_{JX^{\xi}}\epsilon_{\omega}}{2\epsilon_{\omega}})\rbrace
    dvol(\mathfrak m),$$
  where $x\in S$ with $S=Z_{(H)}$ or $S=S_i$ for some $S_i$ as in
  Lemma~\ref{lem-split1} or in Lemma~\ref{lem-split2}, $\xi\in\mathfrak m$, and $\gamma_{\xi}=e^{it\xi}\cdot x$, $t\in [0, 1]$.\\

   \begin{remark}\label{frak m}
   In this and the next subsections, for simplicity, we will only use $\mathfrak m$
   to denote the orthogonal complement of $\mathfrak h$ or of
   $\mathfrak h'$ in $\mathfrak g$ as we did in Formula (\ref{eq1}).
   \end{remark}
  \begin{theorem}\label{estout}
  Consider  $G_{\mathbb C}\cdot S$, where  $S=Z_{(H)}$ or $S=S_i$ for an $S_i$ as in
  Lemma~\ref{lem-split1} or in Lemma~\ref{lem-split2}. There exist constants $b$, and  $D>0$
  such that for all $[x]\in S/G$ (the integral is a function of $[x]$),
  and for all $R$ and  $k$ sufficiently large,

  $$\int_{\mathfrak m-B_R(0)}\tau(\xi, x)exp\lbrace
   -2k\int_{\gamma_{\xi}}\phi_{\xi}\rbrace dvol(\mathfrak m)\leq be^{-RDk},$$
  where $\mathfrak m$ is the orthogonal complement of $\mathfrak h$ or of
   $\mathfrak h'$ in $\mathfrak g$, and $B_R(0)$ is a ball in $\mathfrak m$ of radius $R$ centered at $0$.\\

  Since we can find a uniform bound for $-\frac{\mathcal
    L_{JX^{\xi}}\epsilon_{\omega}}{2\epsilon_{\omega}}$ on $M$, the above inequality
  is also true for the integral
  $\int_{\mathfrak m-B_R(0)}
   \tau(\xi, x)exp\lbrace
   -\int_{\gamma_{\xi}}(2k\phi_{\xi}+\frac{\mathcal
    L_{JX^{\xi}}\epsilon_{\omega}}{2\epsilon_{\omega}})\rbrace
    dvol(\mathfrak m)$.
 \end{theorem}
 The proof of this theorem relies on the following two lemmas.

    \begin{lemma}\label{lem-exp}
   Consider  $G_{\mathbb C}\cdot S$, where  $S=Z_{(H)}$ or $S=S_i$ for an $S_i$ as in
  Lemma~\ref{lem-split1} or in Lemma~\ref{lem-split2}.
  For any $t_0>0$, there exists $C>0$ such that for all $t>t_0$,
  $$exp\lbrace -\int_{\gamma_{t\hat\xi}}2k\phi_{t\hat\xi}\rbrace\leq e^{-2ktC}$$
  uniformly on $S/G$, where $\hat\xi\in\mathfrak m$ with
$|\hat\xi|=1$.
\end{lemma}

 \begin{proof}
  By definition, $\int_{\gamma_{t\hat\xi}}\phi_{t\hat\xi}=\int_0^1<\phi(e^{i\tau t\hat\xi}\cdot x), t\hat\xi>d\tau=t \int_0^1<\phi(e^{i\tau t\hat\xi}\cdot x), \hat\xi>d\tau$.
 Hence, we need to find a positive lower bound for the function  $f_t(\hat\xi, x)=\int_0^1\phi_{\hat\xi}(e^{i\tau t\hat\xi}\cdot x)d\tau$ when $t$ is sufficiently large.
 We prove the lemma for the case $S=Z_{(H)}$. The argument applies
 to other cases.
 Since
 $f_t$ is $G$-invariant,  on each $G$-orbit, we only need to consider a particular point $x$ which has isotropy Lie algebra exactly $\mathfrak h$. So we take
 $S^{\mathfrak h}\subset S$ to be the set of such points. First, fix  $\hat\xi\in\mathfrak m$ with
$|\hat\xi|=1$, and consider the $\hat\xi$-moment map $\phi_{\hat\xi}$.
 Then, $\phi_{\hat\xi}(S^{\mathfrak h})=\mbox{constant}$. For $x\in S^{\mathfrak h}$,
$f_t(\hat\xi, x)$ for any $t>0$ is strictly increasing since
$e^{i\tau t\hat\xi}\cdot x$ is the gradient line of
$\phi_{t\hat\xi}$ and $JX^{t\hat\xi}(x)\neq 0$. If  $S^{\mathfrak
h}$ is compact, we can find a positive lower bound $C_{\hat\xi}$ for
$f_t(\hat\xi, x)$ for all points in  $S^{\mathfrak h}$ and for all
$t>t_0$ for any chosen $t_0>0$. If $S^{\mathfrak h}$ is not compact,
we do the following. Consider a nearby regular
 value  $a>0$ of
$\phi_{\hat\xi}$. For $y\in\phi_{\hat\xi}^{-1}(a)$, consider
 the function   $f_t'(\hat\xi, y)=\int_{-\epsilon}^{1-\epsilon}\phi_{\hat\xi}(e^{i\tau t\hat\xi}\cdot y)d\tau$,
where $\epsilon$ is a small number. By choosing $a$ properly, for
each $x$ in $S^{\mathfrak h}$, there exists
$y\in\phi_{\hat\xi}^{-1}(a)$, such that $f_t'(\hat\xi,
y)=f_t(\hat\xi, x)$ (since the $x$'$s$ are not fixed by the circle
action generated by $\hat\xi$, this can be achieved).
 We choose the positive minimum of $f_t'$ on its
  compact domain $\phi_{\hat\xi}^{-1}(a)$  as $C_{\hat\xi}$ (the positivity of
$C_{\hat\xi}$ is due to $a$ is a regular value). So, for each
$\hat\xi\in\mathfrak m$ with $|\hat\xi|=1$, there exists
$C_{\hat\xi}>0$, such that $f_t(\hat\xi,
x)=\int_0^1\phi_{\hat\xi}(e^{i\tau t\hat\xi}\cdot x)d\tau\geq
C_{\hat\xi}$ for all $[x]\in S/G$. By the compactness of the set
$\{\hat\xi\in \mathfrak m,\, |\hat\xi|=1\}$, and by continuous
dependence of $f_t$ on $\hat\xi$, we can find a positive constant
$C$ such that $f_t\geq C$ uniformly for all $[x]\in S/G$ and for all
$\hat\xi\in\mathfrak m$ with $|\hat\xi|=1$.

   For the proof of other $S$'$s$, we replace the above  $S^{\mathfrak
   h}$ by $S^{\mathfrak h'}\cap\phi^{-1}(a_i)$ (recall that $0\neq
   a_i\in\mathcal F_i$) so that $\phi_{\hat\xi}(S^{\mathfrak
   h'}\cap\phi^{-1}(a_i))=\mbox{constant}$, noticing the fact that
  $S^{\mathfrak h'}\cap\phi^{-1}(a_i)$ has all the representatives
  of $S/G$.
 \end{proof}

 \begin{lemma}
  Consider  $G_{\mathbb C}\cdot S$, where  $S=Z_{(H)}$ or $S=S_i$ for an $S_i$ as in
  Lemma~\ref{lem-split1} or in Lemma~\ref{lem-split2}.
  There exist constants $a$ and $b>0$ such that for all $t>0$
 $$\tau(t\hat\xi, x)\leq bt^{-m}e^{at}$$
 uniformly on $S/G$, where $\hat\xi\in\mathfrak m$ with $|\hat\xi|=1$, and
 $m$ is the dimension of $\mathfrak m$.
\end{lemma}
 \begin{proof}
 The manifold  $G_{\mathbb C}\cdot S$ is a complex
 submanifold of $M$.
 Since $M$ can be embedded into projective spaces,  $G_{\mathbb C}\cdot S$ is a
 complex submanifold of projective spaces.
 The proof of Lemma 5.7 in \cite{HK} applies.
(The proof of  Lemma 5.7 in \cite{HK} does not need the domain of $x$ to be compact, but it  uses the fact that the domain of $\hat\xi$ is compact.)
\end{proof}

  Once we have the above two lemmas, using polar coordinates,  we can prove
  Theorem~\ref{estout}.
  One may refer to the proof of Theorem 5.5 in \cite{HK}.

 \subsection{Approximation}

 \begin{lemma}\label{lem-tau}
  The function $\tau(\xi, x)$ equals $vol(G\cdot x)$ on
  $S$, where  $S=Z_{(H)}$ or $S=S_i$ for an $S_i$ as in
  Lemma~\ref{lem-split1} or in Lemma~\ref{lem-split2}.
 \end{lemma}
  \begin{proof}
  We prove the lemma for the case $S=S_i$ for some $i$. The proof for the other $S$'$s$ is similar.
  Consider the complex submanifold $G_{\mathbb{C}}\cdot S$.
 We take $S^{\mathfrak h'}\subset S$,
 the set of points with isotropy Lie algebra exactly $\mathfrak h'$. Let
 $\tilde S=S^{\mathfrak h'}\cap\phi^{-1}(a_i)$. Then $\tilde S$ contains all the representatives of $S/G$.
        Since $\tau(\xi, x)$ is $G$-invariant,  we only need to
consider the value $\tau(0, x)$ with $x\in \tilde S$.
So we only
consider  Formula
(\ref{eq1}) on $e^{i\mathfrak m}\cdot \tilde S$. 
Consider the submanifold $e^{i\mathfrak m}\cdot\tilde S$. At each point $x\in \tilde S$, the $B$-orthogonal complement of
$T_x \tilde S$ in $T_x (e^{i\mathfrak m}\cdot
  \tilde S)$ is exactly the linear span of the vectors $JX^{\xi}$ with $\xi\in\mathfrak m$: for any $JX^{\xi}$ with $\xi\in\mathfrak m$ and any vector
    $v\in T_x \tilde S$, we have $B(JX^{\xi}, v)_x=\omega(v, X^{\xi})_x=v (\phi_{\xi})_x=0$ since $\phi_{\xi}$ takes
    constant value on  $\tilde S$. So $B$ is block
    diagonalizable at $x$ on the submanifold $e^{i\mathfrak m}\cdot \tilde S$, and
$$dvol(e^{i\mathfrak m} \cdot \tilde S)_x=\sqrt {detB(JX^{\xi_i}, JX^{\xi_j})_x}\,\,dvol(\mathfrak m)\wedge dvol(\tilde S)_x$$
 with $\xi_i, \xi_j\in\mathfrak m$.
  By Lemma~\ref{lem1}, $\sqrt {detB(JX^{\xi_i},
  JX^{\xi_j})_x}=vol(G\cdot x)$.
  \end{proof}
 
 The result of the above lemma will be used in the proof of the following lemma.

  \begin{lemma}\label{limit}
   Consider  $G_{\mathbb C}\cdot S$, where  $S=Z_{(H)}$ or $S=S_i$ for an $S_i$ as in
  Lemma~\ref{lem-split1} or in Lemma~\ref{lem-split2}.
  Define  $${I}_{k, R}([x])=(k/2\pi)^{m/2}\int_{B_R(0)}\tau(\xi,
  x)e^{-kf(\xi, x)}dvol(\mathfrak m),$$
  where $f(\xi, x)=2\int_0^1\phi_{\xi}(e^{it\xi}\cdot x)dt$ at a point $x\in S$ with
  isotropy Lie algebra $\mathfrak h$ or $\mathfrak h'$, $\mathfrak m$ is the orthogonal complement
  of $\mathfrak h$ or of $\mathfrak h'$ in $\mathfrak g$, and $m=$dim$(\mathfrak m)$.

  Then there exists some $R>0$ such that
  $$\mbox{lim}_{k\rightarrow\infty}|{I}_{k,
  R}([x])-2^{-m/2}|=0$$
  uniformly on $S/G$.
   \end{lemma}

  \begin{proof}
  We will prove the lemma for the case $S=Z_{(H)}$. The other cases
 follow similarly.
 We refer to the proof of Lemma 5.10 in \cite{HK}.
 By Theorem~\ref{normvary'},
  the function  $f(\xi, x)$ is a $G$-invariant Morse-Bott function on $G_{\mathbb C}\cdot Z_{(H)}$ with $0\times Z_{(H)}$  being a minimum. By the Morse-Bott lemma, for each point $x\in Z_{(H)}$,
 there exists a neighborhood of this point on which $f(\xi, x)$ can be written
  as a quardratic function. If $Z_{(H)}$ is compact, we can choose the smallest positive
 radius of the (finitely many) neighborhoods as $R$. If $Z_{(H)}$ is not compact,
note that if $Z_{(K)}$ is in the closure of $Z_{(H)}$, then
 (up to conjugacy) the orthogonal complement  $\mathfrak m'$
 of $Lie(K)$  is a linear subspace of the orthogonal complement $\mathfrak m$ of $Lie(H)$.
 Because of this property,
 for the function $f(\xi, x)$ on  $G_{\mathbb C}\cdot Z_{(K)}$,
 we may assume that the neighborhoods of the points $x$'$s\in Z_{(K)}$
  overlap the strata $Z_{(H)}$'$s$
  whose closures contain $Z_{(K)}$. So, we can  use the compactness of $\phi^{-1}(0)$
  to have finitely many neighborhoods, and therefore to
  choose the smallest $R$ for all the strata $Z_{(H)}\subset\phi^{-1}(0)$.
  Once $R$ is chosen, on each $G_{\mathbb C}\cdot Z_{(H)}$, follow
  the arguments of the proof of  Lemma 5.10 in \cite{HK}.  In
 the proof of Lemma 5.10 in \cite{HK}, there are some estimates on the bounds of the absolute value of some continuous functions  of $x\in Z_{(H)}$ which involve certain integrals of the derivative of $\tau(\xi, x)$
 in the direction of $\xi$ (the constants $Q_1$ and $Q_2$). If $Z_{(H)}$ is
 not compact, the formula on $\tau$ in
  (\ref{eq1}) of Section 11.1
  should continuously transform from higher dimensional strata $Z_{(H)}$ of
  $\phi^{-1}(0)$ to lower dimensional ones.
  This should allow us to  extend continuously
  the above continuous functions to the closure of $Z_{(H)}$ in $\phi^{-1}(0)$ and take the maximal
  of the absolute values. (The constant $Q_3$ in the proof of
 Lemma 5.10 in \cite{HK},
 is   $2^{m/2}$ in our case.)
\end{proof}

   \begin{lemma}\label{limit'}
    Consider $G_{\mathbb C}\cdot S$, where $S=Z_{(H)}$ or $S=S_i$ for an $S_i$ as in
  Lemma~\ref{lem-split1} or in Lemma~\ref{lem-split2}.
  Define
  $${J}_{k, R}([x])=(k/2\pi)^{m/2}2^{m/2}\int_{B_R(0)}
   \tau(\xi, x)e^{-kf(\xi, x)}exp\lbrace
   -\int_{\gamma_{\xi}}\frac{\mathcal
    L_{JX^{\xi}}\epsilon_{\omega}}{2\epsilon_{\omega}}\rbrace
    dvol(\mathfrak m).$$
  Then,  there exists  $R>0$ such that
   $$\mbox{lim}_{k\rightarrow\infty}|{J}_{k,R}([x])-1|=0$$
  uniformly on $S/G$.
 \end{lemma}
     \begin{proof}
   In the proof of Lemma~\ref{limit}, replace $\tau(\xi, x)$ by
  $\tau(\xi, x)exp\lbrace
   -\int_{\gamma_{\xi}}\frac{\mathcal
    L_{JX^{\xi}}\epsilon_{\omega}}{2\epsilon_{\omega}}\rbrace$,  just to notice that
 the exponent is $0$ when $\xi=0$.
   \end{proof}

 \subsection{Proof of Theorem~\ref{asym}}

  \begin{proof}
 (a).  We write $I_k^{\mathcal S_{(H)}}$ as the sum of an integral over
$B_R(0)$ and an
  integral over the complement of $B_R(0)$. The result follows from Lemma~\ref{limit}
   and Theorem~\ref{estout}. The proof for $J_k^{\mathcal S_{(H)}}$ is similar but
   using Lemma~\ref{limit'} and Theorem~\ref{estout}.\\

 (b). We assume that $II_k^{Z_{(H)}}\neq 0$ and $\widetilde{II}_k^{Z_{(H)}}\neq
 0$.

  Now we prove $\mbox{lim}_{k\rightarrow\infty}II_k^{Z_{(H)}}=0$.
 Since we have a finite summation in the expression of $II_k^{Z_{(H)}}$,
 we only need to prove that each summand
    goes to $0$ when $k\rightarrow\infty$.
   We will simply write $G_{\mathbb C}\cdot S_i$ as $G_{\mathbb C}\cdot S$.
By Theorem~\ref{estout} and Lemma~\ref{limit}, there exists $K_0>0$,
such
 that when $k>K_0$,
$$\int_{\mathfrak m'}\tau(\zeta, u)exp\lbrace
   -2k\int_{\gamma_{\zeta}}\phi_{\zeta}\rbrace dvol(\mathfrak m')\leq be^{-RDk}+2(k/2\pi)^{-m'/2}2^{-m'/2}$$
 $\leq be^{-RDk}+(k/2\pi)^{-m'/2}b'$.

Now, let us consider the term $\int_{S}|s|^2 (u) dvol(S)$. By
Lemma~\ref{D_2}, we can take $a'\in\mathcal F$ and
 take $S'\subset\phi^{-1}(G\cdot a')$ such that
 $S$ can be reached from $S'$ by following the flow lines  of the vector fields $JX^{\zeta}$, where
 $\zeta\in (\mathfrak m')$.  We use
Theorem~\ref{normvary'} (a) to express $|s|^2 (u)$ in terms of
$|s|^2 (u')$ and we use the arguments in the proof of
Lemma~\ref{lem-exp}  to find a constant $C'>0$ such that $|s|^2
(u)\leq |s|^2 (u')e^{-kC'}$ for all $u'\in S'$. Now, since  $M$ is
compact, $|s|^2 ( u')$ is bounded. The volume of $S$ is also
bounded. So $\int_{S}|s|^2 ( u) dvol(S)\leq C'' e^{-kC'}$ for some
constant $C''$.

 So, for each summand in $II_k^{Z_{(H)}}$, there exist $K_0>0$ and constants
 $C, C', b, b', R, D$  with $C', R$ and  $D$
 positive such that when $k>K_0$, the summand
 $$\leq (k/2\pi)^{n^{(H)}_{\mathfrak{(h')}}/2}Ce^{-kC'}(be^{-RDk}+(k/2\pi)^{-m'/2}b').$$
 Therefore  $\mbox{lim}_{k\rightarrow\infty}II_k^{Z_{(H)}}=0.$

  The proof for the statement about $\widetilde{II}_k^{Z_{(H)}}$ is similar.
 \end{proof}

 \section{Asymptotic unitarity}

  Now, it comes to the  definition of the inner products on $\mathcal H(M, L^{\otimes k})^G$ and
  on $\mathcal H(M, L^{\otimes k}\otimes\sqrt K)^G$. Recall that
   $M^{ss}$ is open and dense in $M$, and
 $$M^{ss}=\bigcup_{(H)} F_{\infty}^{-1} (Z_{(H)}).$$
 There is an open and dense
  set
  $F_{\infty}^{-1} (Z_{(H)})$ for some $H$ in $M^{ss}$, and, correspondingly, there
 is an open and dense stratum $\mathcal S_{(H)}$ in $M//G$. Let us denote the open dense piece  $F_{\infty}^{-1} (Z_{(H)})$ as
 $F_{\infty}^{-1} (Z^O)$, and denote the corresponding open and dense stratum
 $\mathcal S_{(H)}$ of $M//G$ as $\mathcal S^O$.

 \begin{definition}\label{def-qn1}
 Let $s_1, s_2\in\mathcal H(M, L^{\otimes k})^G$ and let $r_1, r_2\in\mathcal H(M, L^{\otimes k}\otimes\sqrt
 K)^G$.
 We define
   $$<s_1, s_2>_{(1)}=\int_M^{(1)}(s_1, s_2) dvol(M)=\int_{F_{\infty}^{-1}
   (Z^O)}(s_1, s_2)
   dvol(F_{\infty}^{-1} (Z^O)),$$
 and  we define
   $$<r_1, r_2>_{(1)}=\int_M^{(1)}(r_1, r_2) dvol(M)=\int_{F_{\infty}^{-1} (Z^O)}(r_1, r_2) dvol(F_{\infty}^{-1} (Z^O)).$$
\end{definition}
By Theorems~\ref{qnsr} and \ref{asym}, we have
\begin{corollary}\label{qnsr1}
Let $s\in\mathcal H(M, L^{\otimes k})^G$, and let  $r\in\mathcal
H(M, L^{\otimes k}\otimes\sqrt K)^G$. Then,
$$\|s\|_{(1)}^2=\int_M^{(1)} |s|^2 dvol(M)=(k/2\pi)^{d_{\mathcal S^O}/2}\int_{\mathcal S^O}|A_k's|^2([x])
I_k^{\mathcal S^O}([x]) \epsilon_{\hat{\omega}_{\mathcal S^O}}+
II_k^{Z^O},$$
 where, $I_k^{\mathcal S^O}([x])=1$ or
 $\mbox{lim}_{k\rightarrow\infty}I_k^{\mathcal S^O}([x])=2^{-d_{G/H}/2}vol(G\cdot
 x)$
 uniformly for $[x]\in \mathcal S^O$  for some $H\neq G$, and,
 $II_k^{Z^O}=0$ or
 $\mbox{lim}_{k\rightarrow\infty}\,II_k^{Z^O}=0$;

$$\|r\|_{(1)}^2=\int_M^{(1)} |r|^2 dvol(M)=(k/2\pi)^{d_{\mathcal S^O}/2}\int_{\mathcal S^O}|B_k'r|^2([x])
J_k^{\mathcal S^O}([x]) \epsilon_{\hat{\omega}_{\mathcal
S^O}}+\widetilde{II}_k^{Z^O},$$
 where, $J_k^{\mathcal S^O}([x])=1$ or $\mbox{lim}_{k\rightarrow\infty}J_k^{\mathcal S^O}([x])=1$ uniformly
 for $[x]\in \mathcal S^O$, and, $\widetilde{II}_k^{Z^O}=0$ or
$\mbox{lim}_{k\rightarrow\infty}\,\widetilde{II}_k^{Z^O}=0$.

\end{corollary}

 The following definition
  modifies the usual definition of quantum norms, but it takes into account all the strata.
 Physical interpretations of this definition would be desirable.

  \begin{definition}\label{def-qn2}
   Let $s_1, s_2\in\mathcal H(M, L^{\otimes k})^G$ and let $r_1, r_2\in\mathcal H(M, L^{\otimes k}\otimes\sqrt
 K)^G$.
  We define
   $$<s_1, s_2>_{(2)}=\int^{(2)}_M (s_1, s_2) dvol(M)=\sum_{Z_{(H)}}\int_{F_{\infty}^{-1}
   (Z_{(H)})}(s_1, s_2)
   dvol(F_{\infty}^{-1} (Z_{(H)})),$$
 and we define
   $$<r_1, r_2>_{(2)}=\int^{(2)}_M (r_1, r_2) dvol(M)=\sum_{Z_{(H)}}\int_{F_{\infty}^{-1} (Z_{(H)})}(r_1, r_2) dvol(F_{\infty}^{-1} (Z_{(H)})).$$
\end{definition}

Again, by Theorems~\ref{qnsr} and \ref{asym}, we have
\begin{corollary}\label{qnsr2}
Let $s\in\mathcal H(M, L^{\otimes k})^G$, and let  $r\in\mathcal
H(M, L^{\otimes k}\otimes\sqrt K)^G$. Then,
$$\|s\|_{(2)}^2=\int^{(2)}_M |s|^2 dvol(M)$$
$$=\sum_{\mathcal S_{(H)}}(k/2\pi)^{d_{\mathcal
S_{(H)}}/2}\int_{\mathcal S_{(H)}}|A_k's|^2([x]) I_k^{\mathcal
S_{(H)}}([x]) \epsilon_{\hat{\omega}_{\mathcal
    S_{(H)}}}+\sum_{Z_{(H)}} II_k^{Z_{(H)}},$$
 where, $I_k^{\mathcal S_{(G)}}=1$ or
  $\mbox{lim}_{k\rightarrow\infty}I_k^{\mathcal S_{(H)}}([x])=2^{-d_{G/H}/2}vol(G\cdot
    x)$ uniformly for $[x]\in \mathcal S_{(H)}$ with $H\neq G$, and,
    $II_k^{Z_{(H)}}=0$ or
    $\mbox{lim}_{k\rightarrow\infty}\,II_k^{Z_{(H)}}=0$;

$$\|r\|_{(2)}^2=\int^{(2)}_M |r|^2 dvol(M)$$
$$=\sum_{\mathcal S_{(H)}}(k/2\pi)^{d_{\mathcal
S_{(H)}}/2}\int_{\mathcal S_{(H)}}|B_k'r|^2([x]) J_k^{\mathcal
S_{(H)}}([x]) \epsilon_{\hat{\omega}_{\mathcal
    S_{(H)}}}+\sum_{Z_{(H)}} \widetilde{II}_k^{Z_{(H)}},$$
 where, $J_k^{\mathcal S_{(G)}}=1$ or $\mbox{lim}_{k\rightarrow\infty}J_k^{\mathcal S_{(H)}}([x])=1$
 uniformly
    for $[x]\in \mathcal S_{(H)}$, and, $\widetilde{II}_k^{Z_{(H)}}=0$ or
    $\mbox{lim}_{k\rightarrow\infty}\, \widetilde{II}_k^{Z_{(H)}}=0$.
\end{corollary}

  For both Definitions~\ref{def-qn1} and \ref{def-qn2}, we have the
  following asymptotic unitarity for the maps $B_k'$.
    \begin{theorem}\label{uni}
   The maps $B_k'$ are asymptotically unitary, in the sense that

   $$\mbox{lim}_{k\rightarrow\infty}\|B_k'^*B_k'-I\|=
   \mbox{lim}_{k\rightarrow\infty}\|B_k'B_k'^*-I\|=0,$$

  where $\|.\|$ refers to the operator norm.
  \end{theorem}
  \begin{proof}
  We use Theorem~\ref{modisom}, the definitions in (\ref{norm_1}) and in (\ref{norm_2}) of Section 5,
   and we use  the above results in Corollaries~\ref{qnsr1} and
  \ref{qnsr2} of Theorems~\ref{qnsr} and \ref{asym}.
  For the case of the quantum norms in
  Definition~\ref{def-qn2}, we also use  the fact that there are finitely many
  strata.
  One may refer to \cite{HK}, the proof of Theorem 5.2 for the asymptotic unitarity of the
  maps
  $B_k$.
  \end{proof}

\end{document}